\newif\ifarxiv
\arxivtrue

\ifarxiv

\documentclass[12pt]{amsart}
\usepackage{amsmath,amssymb}
\usepackage{mathptmx}
\usepackage{graphicx}
\usepackage{amscd}
\usepackage{overpic}
\usepackage{fullpage}
\else 

\documentclass[11pt]{article}
\usepackage{fullpage}

\usepackage{amsmath,amssymb, amsthm, amscd}
\usepackage{graphicx}
\usepackage{overpic}

\fi 

\def\RCS$#1: #2 ${\expandafter\def\csname RCS#1\endcsname{#2}}
\RCS$Revision: 1.82 $
\RCS$Date: 2005/12/07 17:01:10 $
\ifarxiv

\else
\relax
\fi

\newcommand{\C}{{\mathbb C}}
\newcommand{\Q}{{\mathbb Q}}
\newcommand{\Z}{{\mathbb Z}}
\newcommand{\maps}{\colon\thinspace}

\newcommand{\abs}[1]{{\left| #1 \right|}}
\renewcommand\b{\beta}
\newcommand\cp{{\bf \mathbb{C}P}}
\newcommand\Sp{{\bf S}}

\newcommand\CF{{\mathcal F}}
\newcommand\CL{{\mathcal L}}
\newcommand\CH{{\mathcal H}}
\newcommand\CM{{\mathcal M}}
\newcommand\CO{{\mathcal O}}
\newcommand\CP{{\mathcal P}}

\newcommand{\hk}{\mathcal{H}}
\newcommand{\khr}{\mathit{KhR}}
\newcommand{\kh}{\overline{\mathit{KhR}}}
\newcommand{\khrn}{\khr_N}
\newcommand{\khs}{\khr_{2,s}}
\newcommand{\Pred}{\CP}
\newcommand{\Pun}{\bar{\CP}}
\newcommand{\pred}{P}
\newcommand{\pun}{\bar{P}}
\newcommand{\HFK}{\widehat{\mathit{HFK}}}
\newcommand{\hfhat}{\widehat{\mathit{HF}}}
\newcommand{\hfk}{\mathit{HFK}}
\newcommand{\hkr}{\mathit{HKR}}
\newcommand{\hkun}{\overline{\mathit{HKR}}}

\newcommand{\dotdia}[4]{\begin{center} 
\includegraphics[scale=#4]{#1}

\vspace{0.5cm}

{\Large $#2 \maps a_{\min} = #3$}
\end{center}}

\newcommand{\dotdianod}[4]{\begin{center} 
\includegraphics[scale=#4]{#1}

\vspace{0.5cm}

{\Large $#2 \maps a_{\min} = #3$}

\vspace{0.1cm}

($d_1$ and $d_{-1}$ are not shown.)

\end{center}}

\ifarxiv
\makeatletter
\@ifundefined{swappedhead@plain}{%
\def\swappedhead#1#2#3{%
   \thmnumber{\@upn{\@secnumfont#2\@ifnotempty{#1}{.~}}}%
   \thmname{#1}%
    \thmnote{ {\the\thm@notefont(#3)}}}
\let\swappedhead@plain=\swappedhead}
\makeatother

\else
\relax
\fi

\theoremstyle{plain} 

\swapnumbers
\newtheorem{theorem}{Theorem}[section]
\newtheorem{conjecture}[theorem]{Conjecture}
\newtheorem{conj}[theorem]{Conjecture}
\newtheorem{lemma}[theorem]{Lemma}

\newtheorem{prop}[theorem]{Proposition}

\newtheorem{goal}[theorem]{Goal}

\theoremstyle{definition}
\newtheorem{definition}[theorem]{Definition}

\theoremstyle{remark}
\newtheorem{remark}[theorem]{Remark}

\newtheorem{example:unknot}[theorem]{Example: The Unknot}
\newtheorem{example:T23}[theorem]{Example: The Trefoil}
\newtheorem{example:T34}[theorem]{Example: \(T_{3,4}\)}
\newtheorem{example:figure8}[theorem]{Example: The Figure-eight Knot}

\makeatletter
  \let\c@theorem=\c@subsection
  \let\c@figure=\c@subsection
  \let\p@figure=\p@subsection
  
  \let\cl@figure=\cl@subsection
  \let\c@table=\c@subsection
  \let\p@table=\p@subsection
  
  \let\cl@table=\cl@subsection
\makeatother

\def\mathcenter#1{%
  \vcenter{\hbox{#1}}%
}

\newcommand{\mfig}[2][]{
        \mathcenter{\includegraphics[#1]{#2}}
}

\ifarxiv 
\begin{document}

\enlargethispage{1.5cm}

\title{The Superpolynomial for Knot Homologies}

\author[Dunfield]{Nathan M.~Dunfield}
\thanks{Mathematics 253-37, California Institute of Technology, Pasadena, CA 91125, USA.  \texttt{dunfield@caltech.edu}}

\author[Gukov]{Sergei Gukov} \thanks{
Jefferson Physical Laboratory, Harvard Univ., Cambridge, MA, 02138, USA.  \texttt{gukov@schwinger.harvard.edu}}

\author[Rasmussen]{Jacob Rasmussen} \thanks{Dept.~of Mathematics,
  Princeton University, Princeton, NJ 08544, USA.  \texttt{jrasmus@math.princeton.edu}}


\begin{abstract}
  We propose a framework for unifying the $sl(N)$ Khovanov-Rozansky
  homology (for all $N$) with the knot Floer homology.  We argue that
  this unification should be accomplished by a triply graded homology
  theory which categorifies the HOMFLY polynomial.  Moreover, this
  theory should have an additional formal structure of a family of
  differentials.  Roughly speaking, the triply graded theory by itself
  captures the large $N$ behavior of the $sl(N)$ homology, and
  differentials capture non-stable behavior for small $N$, including
  knot Floer homology.  The differentials themselves should come from
  another variant of $sl(N)$ homology, namely the deformations of it
  studied by Gornik, building on work of Lee.  
  
  While we do not give a mathematical definition of the triply graded theory,
  the rich formal structure we propose is powerful enough to make many
  non-trivial predictions about the existing knot homologies that can
  then be checked directly.  We include many examples where we can
  exhibit a likely candidate for the triply graded theory, and these
  demonstrate the internal consistency of our axioms.  We conclude
  with a detailed study of torus knots, developing a picture which
  gives new predictions even for the original $sl(2)$ Khovanov
  homology.
\end{abstract}

\maketitle
\thispagestyle{empty}

\else
\title{The Superpolynomial for Knot Homologies}

\author{Nathan M.~Dunfield, Sergei Gukov, and Jacob Rasmussen \\
Mathematics 253-37, California Institute of Technology, Pasadena, CA 91125, USA. \\
Jefferson Physical Laboratory, Harvard Univ., Cambridge, MA, 02138, USA. \\ 
Dept.~of Mathematics, Princeton University, Princeton, NJ 08544, USA.  \\
\texttt{dunfield@caltech.edu, gukov@schwinger.harvard.edu, jrasmus@math.princeton.edu}
}

\begin{document}

\maketitle 
\begin{abstract}
  We propose a framework for unifying the $sl(N)$ Khovanov-Rozansky
  homology (for all $N$) with the knot Floer homology.  We argue that
  this unification should be accomplished by a triply graded homology
  theory which categorifies the HOMFLY polynomial.  Moreover, this
  theory should have an additional formal structure of a family of
  differentials.  Roughly speaking, the triply graded theory by itself
  captures the large $N$ behavior of the $sl(N)$ homology, and
  differentials capture non-stable behavior for small $N$, including
  knot Floer homology.  The differentials themselves should come from
  another variant of $sl(N)$ homology, namely the deformations of it
  studied by Gornik, building on work of Lee.  
  
  While we do not give a mathematical definition of the triply graded theory,
  the rich formal structure we propose is powerful enough to make many
  non-trivial predictions about the existing knot homologies that can
  then be checked directly.  We include many examples where we can
  exhibit a likely candidate for the triply graded theory, and these
  demonstrate the internal consistency of our axioms.  We conclude
  with a detailed study of torus knots, developing a picture which
  gives new predictions even for the original $sl(2)$ Khovanov
  homology.
\end{abstract}

\fi

\setcounter{tocdepth}{1}
\tableofcontents


\section{Introduction}

\subsection{Knot homologies}

Here, we are interested in homology theories of knots in $\Sp^3$
associated to the HOMFLY polynomial.  For a knot $K$, its HOMFLY
polynomial $\pun(K)$ is determined by the skein relation:
\begin{equation}
a \pun \left( \, \mfig[scale=0.6]{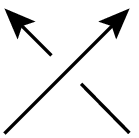} \, \right)
- a^{-1} \pun  \left( \, \mfig[scale=0.6]{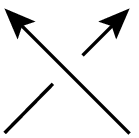} \,  \right)
= ( q - q^{-1}) \pun \left( \, \mfig[scale=0.6]{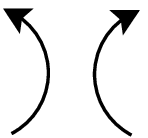} \, \right),
\end{equation}
together with the requirement that $\pun(\mathrm{unknot}) = (a -
a^{-1})/(q-q^{-1})$.  The HOMFLY polynomial unifies the quantum
$sl(N)$ polynomial invariants of $K$, which are denoted by
$\pun_N(K)(q)$ and are equal to $\pun(K)(a = q^N, q)$.
Here, the original Jones polynomial $J(K)$ is just $\pun_2(K)$.
The HOMFLY polynomial encodes the classical Alexander
polynomial as well.

A number of different knot homology theories have been discovered
related to these polynomial invariants.  Although the details of these
theories differ, the basic idea is that for a knot $K$, one can
construct a doubly graded homology theory $H_{i,j} (K) $ whose graded
Euler characteristic with respect to one of the gradings gives a
particular knot polynomial.  Such a theory is referred to as a
\emph{categorification} of the knot polynomial.

For example, the Jones polynomial $J$ is the graded Euler
characteristic of the doubly graded \emph{Khovanov Homology}
$H_{i,j}^{\mathit{Kh}}(K)$; that is,
\begin{equation}
J (q) = \sum_{i,j} (-1)^j q^i \dim H_{i,j}^{\mathit{Kh}} (K).
\end{equation}
Here, the grading $i$ is called the \emph{Jones grading}, and $j$ is called the \emph{homological grading}.    Khovanov originally
constructed $H_{i,j}^{\mathit{Kh}}$ combinatorially in terms of skein
theory \cite{Khovanov}, 
but it is conjectured to be essentially the same as Seidel and
Smith's \emph{symplectic Khovanov homology} which is defined by
considering the Floer homology of a certain pair of Lagrangians
\cite{SeidelSmith}.

Khovanov's theory was generalized by Khovanov and Rozansky
\cite{RKhovanov} to categorify the quantum $sl(N)$ polynomial
invariant $\pun_N(q)$.  Their homology $\hkun_{i,j}^N(K)$
satisfies
\begin{equation}
\pun_N (q) = \sum_{i,j} (-1)^j q^i \dim \hkun_{i,j}^{N} (K).
\end{equation}
For $N=2$, this theory is expected to be equivalent
to the original Khovanov homology.
There are also important deformations of the original Khovanov
homology \cite{ESL2,DBN04,KhovanovFrob},
as well as of the $sl(N)$ Khovanov-Rozansky homology \cite{Gornik}.
In a sense, the deformed theory of Lee \cite{ESL2} also can be regarded
as a categorification of the $sl(1)$ polynomial invariant.

Another knot homology theory that will play an important role here is
knot Floer homology, $\HFK_j (K; i)$, introduced in \cite{OShfk,Rasmussen}.
It provides a categorification of the Alexander polynomial:
\begin{equation}
\Delta (q) = \sum_{i,j} (-1)^j q^i \dim \HFK_j (K;i).
\end{equation}
Unlike Khovanov-Rozansky homology, knot Floer homology is not
known to admit a combinatorial definition; in the end, computing
$\HFK$ involves counting pseudo-holomorphic curves.

The polynomials above are closely related as they can all be derived
from a single invariant, namely the HOMFLY polynomial.  While the
above homology theories categorify polynomial knot invariants in the
same class, their constructions are very different!
Despite this, our objective here is
\begin{goal}\label{main-goal}
Unify the Khovanov-Rozansky $sl(N)$ homology (for all $N$),
knot Floer homology, and various deformations thereof
into a single theory.
\end{goal}
We do not succeed here in defining such a unified theory.
Instead, we postulate a very detailed picture of what such
a theory should look like: it is a triply graded homology
theory categorifying the HOMFLY polynomial together with
a certain additional formal structure.
Although we don't know a definition of this triply graded theory,
our description of its properties is powerful enough to give us many
non-trivial predictions about knot homologies that can be verified directly.


\begin{figure}
\begin{center}
\includegraphics[width=3.5in]{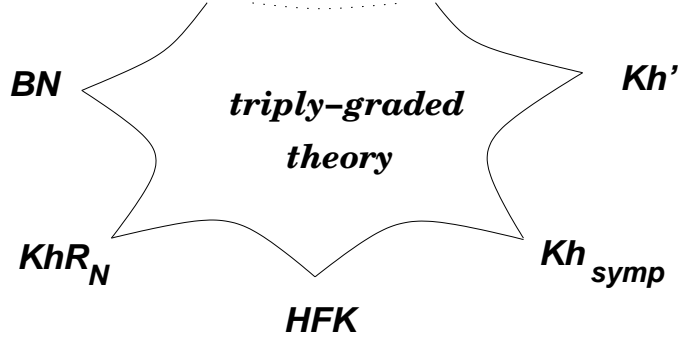}
\end{center}
\caption{Triply-graded theory as a unification of knot homologies.}
\label{Fig:superpol}
\end{figure}

There are several reasons to hope for the type of unified theory asked
for in Goal~\ref{main-goal}. In the recent work \cite{GSV},
a physical interpretation of the Khovanov-Rozansky homology
naturally led to the unification of the $sl(N)$ homologies,
when $N$ is sufficiently large. At the small $N$ end, the $sl(2)$
Khovanov homology and $\HFK$ seem to be very closely related.
For instance, their total ranks are very often (but not always)
equal (see \cite{KFC} for more).
One hope for our proposed theory is that it will explain
the mysterious fact that while the connections between
$\hkun_2$ and $HFK$ hold very frequently, they are not universal.


\subsection{The superpolynomial}\label{subsec-intro-super}

We now work toward a more precise statement of our proposed
unification, starting with a review of the work \cite{GSV}.
To concisely describe the homology groups
$\hkun_{i,j}^N(K)$, it will be useful to introduce the graded Poincar\'e
polynomial, $\kh_N (q,t) \in \Z [q^{\pm 1}, t^{\pm 1}]$, which encodes the
dimensions of these groups via
\begin{equation}
\label{khndef}
\kh_N (q,t) := \sum_{i,j} q^i t^j  \dim \hkun_{i,j}^{N} (K).
\end{equation}
The Khovanov-Rozansky homology has finite total dimension, so $\kh_N$
is a \emph{finite polynomial}, that is, one with only finitely many non-zero
terms.  The Euler characteristic condition on $\hkun_{i,j}^N(K)$ is
concisely expressed by $\pun_N(q) = \kh_N(q, t = -1)$.

The basic conjecture of \cite{GSV} is that
\begin{conj}\label{khnconj}
  There exists a finite polynomial $\Pun (K) \in \Z [a^{\pm
    1},q^{\pm 1},t^{\pm 1}]$ such that
\begin{equation}
\label{khnviap}
\kh_N (q,t) = \frac{1}{q-q^{-1}} \Pun (a=q^N,q,t)
\end{equation}
for all sufficiently large $N$.
\end{conj}

We will refer to $\Pun(K)$ as the \emph{superpolynomial} for $K$.
This conjecture essentially says that, for sufficiently large $N$, the
dimension of the $sl(N)$ knot homology grows linearly in $N$, and the
precise form of this growth can be encoded in a finite set of the
integer coefficients.  Therefore, if one knows the $sl(N)$ knot
homology for two different values of $N$, both of which are in the
``stable range'' $N \geq N_0$, one can use \eqref{khnviap} to determine
the $sl(N)$ knot homology for all other values of $N \geq N_0$.

In some examples, it seems that \eqref{khnviap} holds true for
all values of $N$, not just large $N$.  In \cite{GSV}, this was used
to compute $\Pun(K)$ for certain knots.  However, this is not always
true.  The simplest knot for which \eqref{khnviap} holds for all $N \geq
3$ but not for $N=2$ is the 8-crossing knot $8_{19}$.
Notice, the Conjecture~\ref{khnconj} implies that, for all knots,
the HOMFLY polynomial is a specialization of the superpolynomial,
\begin{equation}
\label{homflyviap}
\pun (K)(a,q) = \frac{1}{q-q^{-1}} \Pun (a,q,t=-1).
\end{equation}

The motivation for Conjecture~\ref{khnconj} in \cite{GSV}
was based on the geometric interpretation of the $sl(N)$ knot
homology and the 3-variable polynomial $\Pun(a,q,t)$.
In fact, we can offer two (related) geometric interpretations
of $\Pun(a,q,t)$:

\begin{itemize}

\item as an index ({\it cf.} elliptic genus):
\begin{equation}
\Pun (a,q,t) = {\rm Str}_{\CH} [a^Q q^s t^r]
= {\rm Tr}_{\CH} [(-1)^F a^Q q^s t^r].
\end{equation}
Here $\CH = \CH_{BPS}$ is a
$\Z_2 \oplus \Z \oplus \Z \oplus \Z$-graded Hilbert space
of the so-called BPS states. Specifically, $F$ is the $\Z_2$ grading,
and $Q$, $s$, and $r$ are the three $\Z$ gradings.
Following the notations in \cite{GSV},
we also introduce the graded dimension
of this Hilbert space:
\begin{equation}\label{dqsrdef}
D_{Q,s,r} := (-1)^F \dim \CH_{BPS}^{F,Q,s,r}.
\end{equation}
Notice that the integer coefficients of the polynomial
$\Pun (a,q,t)$ are precisely the graded dimensions \eqref{dqsrdef}.

\item as an enumerative invariant: the triply graded
integers $D_{Q,s,r}$ are related to the dimensions of
the cohomology groups:
\begin{equation}\label{hmgq}
H^k (\CM_{g,Q})
\end{equation}
where $\CM_{g,Q}$ is the moduli space of holomorphic Riemann surfaces
with boundary in a certain Calabi-Yau 3-fold.  We will come back to
this relationship in Section~\ref{sec-geometry}.

\end{itemize}


\subsection{Reduced superpolynomial}\label{sec-intro-redsuper}

The setup of the last section needs to be modified in order to bring
knot Floer homology into the picture.  Let $\pred(K)(a,q)$ be the
\emph{reduced} or \emph{normalized} HOMFLY polynomial of the knot $K$,
determined by the convention that $\pred(\mathrm{unknot})=1$.  This switch
brings the Alexander polynomial naturally into the picture since it
arises by a specialization $\Delta(q) = \pred(K)(a = 1, q)$.  There is a
categorification of $\pred(K)(a = q^N, q)$ called the \emph{reduced}
Khovanov-Rozansky Homology (see \cite[\S3]{Kpatterns} and
\cite[\S7]{RKhovanov}).  We will use $\khr_N(K)(q,t)$ to denote the
Poincar\'e polynomial of this theory.

For this reduced theory, there is also a version of the Conjecture
\ref{khnconj}.  Essentially, it says that, for sufficiently large $N$,
the total dimension of the reduced $sl(N)$ knot homology is
independent of $N$, and the graded dimensions of the homology groups
change linearly with $N$:
\begin{conj}
\label{khredconj}
There exists a finite polynomial $\Pred (K) \in \Z_{\geq 0} [a^{\pm 1},q^{\pm
  1},t^{\pm 1}]$ such that
\begin{equation}
\label{khviapred}
\khr_N (q,t) = \Pred (a=q^N,q,t)
\end{equation}
for all sufficiently large $N$.
\end{conj}

In contrast with the previous case, in the reduced case the
superpolynomial is required to have non-negative coefficients.  This
is forced merely by the form of \eqref{khviapred}, since for large $N$
distinct terms in $\Pred(a,q,t)$ can't coalesce when we specialize to
$a = q^N$.  Moreover, one also has
\begin{equation}
\label{homflyviapred}
\pred(K)(a,q) = \Pred(a, q, t=-1).
\end{equation}
Thus we will view $\Pred(a,q,t)$ as the Poincar\'e polynomial of some
new triply graded homology theory $\hk_{i,j,k}(K)$ categorifying 
the normalized HOMFLY polynomial.

As with unreduced theory, for some simple cases \eqref{khviapred}
holds for all $N \geq 2$.  However, in general there will be exceptional
values of $N$ for which this is not the case.  To account for this, we
introduce an additional structure on $\hk_*(K)$, a family of
differentials $\{ d_N \}$ for $N > 0$.  The complete details of this
structure we postpone until Section~\ref{sec-differentials}, but the
basic idea is this.  The $sl(N)$ homology is the homology of
$\hk_*(K)$ with respect to the differential $d_N$.  For large $N$, the
differential $d_N$ is trivial, giving the stabilization phenomena of
Conjecture~\ref{khredconj}.  The main reason for expecting the
presence of the differentials \(d_N\) for comes from Gornik's work on
the \(sl(N)\) homology.  In particular, in \cite{Gornik} Gornik
describes a deformation of Khovanov and Rozansky's construction which
gives rise to a differential on \(\hkr_N\).

We also postulate additional differentials for $N \leq 0$.
After a somewhat mysterious re-grading, the knot Floer homology
arises from the $N=0$ differential. Consider the Poincar\'e polynomial
\begin{equation}
\label{hfkdef}
\hfk (q,t) := \sum_{i,j} q^i t^j  \dim \HFK_j (K; i).
\end{equation}
In the simplest cases, we have the following relationship between
the knot Floer homology and the superpolynomial:
\begin{equation}
  \hfk(q, t) = \Pred(a = t^{-1}, q, t).
\end{equation}
For the more general situation, see Section~\ref{sec-differentials}.  

\subsection{Some Predictions}

Our conjectures imply that the HOMFLY polynomial, the knot Floer
homology, and Khovanov-Rozansky homology should all be related.
Unfortunately, this relation is mediated by the triply-graded homology
group \(\hk_{i,j,k}(K)\), which is often considerably larger than
\(\HFK(K)\), \(\hkr_2(K)\), or the minimum size dictated by \(P(K)\).
Thus it seems unlikely that there will be a general relation between
the dimensions of either of these groups and the HOMFLY polynomial.
On the other hand, our hypotheses about the structure of the triply
graded theory enable us to make testable predictions about the $sl(2)$
Khovanov homology and HOMFLY polynomial for some specific families of
knots. We list some of the more important ones here:

\begin{enumerate}
\item {\bf \(\hkr_N\) for small knots:} 
Using conjectured properties of the
  triply graded theory, we make exact predictions for the
  the group \(\hk(K)\) for many knots with 10 crossings or
  fewer. These are given in Sections~\ref{Sec:Examples} and
  \ref{sec-dot-diagrams}.
From them, it is easy to predict the form of \(\khrn(K)\)
  for \(N>2\). These predictions have been verified in simple
  cases \cite{JakeUnwritten}; to check them in  others requires better
  methods for calculating the Khovanov-Rozansky homology. 

\item {\bf HOMFLY polynomials of thin knots:} In
  Section~\ref{SubSec:Thin}, we describe a class of {\it \(\hk\)-thin} knots
  whose triply-graded homology has an especially simple form. Let \(K\)
  be such a knot, and let \(T\) be the \((2,n)\) torus knot
  with the same signature as \(K\). Then our conjectures imply that 
the quotient
\begin{equation*}
\frac{P(K) -  P(T)}{(1-a^2q^2)(1-a^2q^{-2})}
\end{equation*}
should be an alternating polynomial. Two-bridge knots are expected to
be \(\hk\)-thin; we have verified that the relation above holds for all
such knots with determinant less than \(200\).

\item { \bf A new pairing on Khovanov homology:}
Our conjectures suggest that for many knots, the
Khovanov polynomial should have the following form:
\begin{equation*}
\khr_2(K) = q^mt^n + (1+q^6t^3)Q_-(q,t)
\end{equation*} 
 where \(Q_-\) is a polynomial with positive coefficients. (See
 Section~\ref{Subsec:d-1} for a complete discussion.) This pattern is
 easily verified in  examples, but so far as we are aware, it had
 previously gone unnoticed. 

\item {\bf Khovanov homology of torus knots:} In
  Sections~\ref{sec-torus-1} and \ref{sec-torus-2}, we use our
  conjectures to make predictions about the \(N=2\) Khovanov homology
  of torus knot which can be checked
  against  the computations made by Bar-Natan \cite{DBNWeb}. These
  predictions provide some of the the best evidence in favor of the
  conjectures, since the Khovanov homology of torus knots had
  previously seemed quite mysterious. 

\end{enumerate}

\subsection{Candidate theories for the superpolynomial}

The most immediate question raised by Conjecture~\ref{khredconj} is
how to define the underlying knot homology whose Poincar{\'e} polynomial
is the the superpolynomial.  In formulating our conjectures, the
approach we had in mind was simply to take the inverse limit of
\(\khrn\) as \(N \to \infty\). This method rests on two basic principles.
First, we should have some sort of map from the \(sl(N)\) homology to
the \(sl(M)\) homology for \(M<N\), and second, for a fixed knot \(K\)
the dimension of \(\hkr_N(K)\) should be bounded independent of \(N\).
We expect that the maps required by the first principle should be
defined using the work of Gornik \cite{Gornik}, although at the
moment, technical difficulties prevent us from giving a complete proof
of their existence. The proof of the second principle should be more
elementary --- it should be essentially skein theoretic in nature.

Very recently, Khovanov and Rozansky have introduced a triply
graded theory categorifying the HOMFLY polynomial \cite{RKhovanovII},
which gives another  candidate for our proposed theory.
This theory has some obvious advantages over the approach
described above; it is already known to be well-defined,
and its definition is in many respects simpler than that
of the \(sl(N)\) theory. At the same time, there are some
gaps between what the theory provides and what our conjectures
suggest that it should have. The most important of these is
the family of differentials \(d_N\) alluded to above.
One of our aims in writing this paper is to encourage
people to look for these differentials, and, with luck, to find them!

Another approach to constructing a knot homology associated
the superpolynomial might be
based on an algebraic structure which unifies $sl(N)$ (or $gl(N)$)
Lie algebras (for all $N$).
A natural candidate for such structure is the infinite dimensional
Lie algebra, $gl(\lambda)$, introduced by Feigin \cite{Feigin}
as a generalization of $gl(N)$ to non-integer, complex values
of the rank $N$. It is defined as a Lie algebra of the following
quotient of the universal enveloping algebra of $sl(2)$:
\begin{equation}
gl(\lambda) = \mathit{Lie} \left(
U(sl(2)) \Big/ C - \frac{\lambda (\lambda-1)}{2} \right)
\end{equation}
where $C$ is the Casimir operator in $U(sl(2))$.
One can also define $gl(\lambda)$ as a Lie algebra
of differential operators on $\cp(1)$
of ``degree of homogeneity" $\lambda$:
\begin{equation}
gl(\lambda) = \mathit{Lie} (\mathit{Diff}_{\lambda})
\end{equation}
Representation theory of $gl(\lambda)$ is very simple and has
all the properties that we need: For generic $\lambda \in \C$,
$gl(\lambda)$ has infinite dimensional representations.
Characters of these representations appear in the superpolynomial
of torus knots! On the other hand, for $\lambda = N$,
we get the usual finite dimensional representations of $gl(N)$.

\subsection{Generalizations}

We expect many generalizations of this story. Thus, from
the physics point of view, it is clear that a categorification
of the quantum $sl(N)$ invariant should exist for arbitrary
representation of $U_q (sl(N))$, not just the fundamental
representation.

\subsection{Contents of the paper}

In the next section we summarize our conventions and notations.  In
Section \ref{sec-differentials}, we introduce families of graded
differentials, which play a key role in the reduction to different
knot homologies, and give a precise statement of our main conjecture.
In Section \ref{sec-geometry}, we explain the geometric interpretation
of the triply graded theory.  Various examples and patterns are
discussed in Section \ref{Sec:Examples}; these serve to illustrate the
internal consistency of our proposed axioms.
Section~\ref{sec-torus-1} begins our study of torus knots, and there
we give a complete conjecture for the superpolynomials of $(2,n)$ and
$(3,n)$ torus knots.  While we don't have a complete picture for
general $(n,m)$ torus knots, in Section~\ref{sec-torus-2} we suggest a
limiting ``stable'' picture as $m \to \infty$.  Finally,
Section~\ref{sec-dot-diagrams} gives information about the
superpolynomial for certain $10$ crossing knots discussed
in Section~\ref{Sec:Examples}.  

\subsection*{Acknowledgments}  

We are grateful to P.~Etingof, B.~Gornik, V.~Kac, M.~Khovanov,
C.~Manolescu, P. Ozsv{\'a}th, A.~Schwarz,
C.~Taubes, C.~Vafa, and Z.~Szab{\'o} for valuable discussions.  N.D. was
partially supported by NSF grant \#DMS-0405491 and a Sloan Fellowship.
This work was conducted during the period S.G. served as a Clay
Mathematics Institute Long-Term Prize Fellow.  J.R. was partially
supported by an NSF Postdoctoral Fellowship.


\section{Notation and conventions}\label{sec-conventions}

In this section, we give our conventions for knot polynomials, and the
various homology theories.  Some of these differ from standard
sources; in particular, we view the $sl(N)$ theory as homology rather
than cohomology. Also, our convention for the knot Floer homology
is the mirror of the standard one.
  The notation used throughout the paper is collected
in Table~\ref{table-notation}.

\newcommand{\mel}{\\[6pt]}
\newcommand{\showcom}[1]{}

\begin{table}
\begin{tabular}{c cp{4.9in}}

$\pred(K)(a, q)$  &  \showcom{pred} & The normalized HOMFLY polynomial of the knot $K$, where $P(\mathrm{unknot}) = 1$. \mel

$\pun(K)(a, q)$  &  \showcom{pun} & The unnormalized HOMFLY polynomial of the knot $K$, where $\pun(\mathrm{unknot}) = (a - a^{-1})/(q - q^{-1})$. \mel

$\hkr^N_{i,j}(K)$ & \showcom{hkr} & The \emph{reduced} $sl(N)$ Khovanov-Rozansky homology of the knot $K$ categorifying $P(K)$. Here $i$ is the $q$-grading and $j$ the homological grading.   \mel

$\hkun^N_{i,j}(K)$ & \showcom{hkun} & The \emph{unreduced} $sl(N)$ Khovanov-Rozansky homology of the knot $K$ categorifying $P(K)$. Here $i$ is the $q$-grading and $j$ the homological grading.  \mel

$\khr_N(K)(q, t)$ & \showcom{khr} & The Poincar\'e polynomial of the \emph{reduced} $sl(N)$ Khovanov-Rozansky homology of the knot $K$.   In particular, $\khr_N(q, t=-1) = P(a=q^N, q)$.   \mel

$\kh_N(K)(q,t)$ & \showcom{kh} & The Poincar\'e polynomial of the \emph{unreduced} $sl(N)$ Khovanov-Rozansky homology of the knot $K$.   In particular, $\kh_N(q, t=-1) = \pun(a=q^N, q)$.   \mel

$\hk_{i,j,k}(K)$ & \showcom{hk} & A triply graded homology theory which categorifies $\pred(K)$.   The indices $i$ and $j$ correspond to the variables $a$ and $q$ of $\pred(K)$ respectively, and $k$ is the homological grading.  \mel 

$\Pred(K)(a,q,t)$ & \showcom{Pred} & The Poincar\'e polynomial of $\hk_{*}(K)$, called the reduced superpolynomial of $K$.  In particular, $\Pred(K)(a,q, t= -1) = \pred(a,q)$.  \mel

$\Pun(K)(a,q,t)$ & \showcom{Pun} & The \emph{unreduced} superpolynomial of the knot $K$.  This is the Poincar\'e polynomial of a triply graded theory  categorifying $\pun(K)$. \mel 

$\Pred_N(q,t)$ & \showcom{Pred} & The Poincar\'e polynomial of the homology of $\hk_{*}(K)$ with respect to the differential $d_N$.  \mel

$\Delta(K)(q)$ &  \showcom{Delta} & The Alexander polynomial of the knot $K$.  With our conventions, it is a polynomial in $q^2$ and is equal to $P(a=1, q)$.   \mel 

$\HFK(K)$ & \showcom{HKF} & The knot Floer homology of the knot $K$.  \mel

$\hfk(K)(q, t)$ & \showcom{hfk} & The Poincar\'e polynomial of $\HFK(K)$, with $q$ corresponding to the Alexander grading, and $t$ the homological grading.  \mel
\end{tabular}
\caption{Notation.}\label{table-notation}
\end{table}

\subsection{Crossings}  Our conventions for crossings are given below:
\[
\mbox{positive} \  =  \mfig[scale=0.6]{L+.eps}  \quad \quad \quad \quad \quad  \mbox{negative} \   =  \mfig[scale=0.6]{L-.eps}
\]
This convention agrees with \cite{GSV}, but differs from \cite[Fig.~8]{Khovanoviii} and \cite[Fig.~45]{RKhovanov}.  

\subsection{Torus knots}  

The torus knot $T_{a,b}$ is the knot lying on a standard solid torus
which wraps $a$ times around in the longitudinal direction and $b$
times in the meridian direction.  For us, the standard $T_{a,b}$ has
\emph{negative} crossings.  In particular, the trefoil knot $3_1$ in
the standard tables \cite{Rolfsen, KnotAtlas} is exactly $T_{2,3}$
with our conventions.  However, it is important to note that some
other torus knots in these tables are positive rather than negative
(\emph{e.g.}~$8_{19}$ and $10_{124}$), and this is why the
superpolynomial for $10_{124}$ given in Section~\ref{sec-dot-diagrams}
differs from Section~\ref{sec-torus-1}.  

\subsection{Signature} 

Our choice of sign for the signature $\sigma(K)$ of a knot $K$ is such
that the $\sigma( T_{2,3}) = 2$.  That is, negative knots have positive
signatures.

\subsection{Knot polynomials}  

For us, the normalized HOMFLY polynomial $\pred$ of an oriented link $L$
is determined by the skein relation
\begin{equation}
a \pred \left( \, \mfig[scale=0.6]{L+.eps} \, \right)  - a^{-1} \pred  \left( \, \mfig[scale=0.6]{L-.eps} \,  \right)  = ( q - q^{-1}) \pred \left( \, \mfig[scale=0.6]{L0.eps} \, \right),
\end{equation}
together with the requirement that $\pred ( \mathrm{unknot} ) = 1$.  The
unnormalized HOMFLY polynomial $\pun(L)$ is determined by the
alternate requirement that $\pun( \mathrm{unknot}) = (a - a^{-1})/(q -
q^{-1})$.

Several different conventions for the HOMFLY polynomial can be found
in the literature; another common one involves the
change $ a \to a^{1/2}, q \to q^{1/2}$.  Also, sources sometimes
simultaneously switch $a \to a^{-1}$ and $q \to q^{-1}$.  For the
negative torus knot $T_{2,3}$, the polynomial $\pred(T_{2,3})$ has all
positive exponents of $a$.

For knots, our conventions are consistent with \cite{GSV} (for links,
the skein relation here differs by a sign).  The papers of Khovanov
and Rozansky \cite{Khovanov, Kpatterns, RKhovanov,
  RKhovanovII} use the convention where $a$ and $q$ are replaced with
their inverses.  For the Knot Atlas \cite{KnotAtlas}, the conventions
for HOMFLY agree with ours if you substitute $z = q - q^{-1}$;
however, the Knot Atlas' conventions for the Jones polynomial differ
from ours by $q \to q^{-1}$.

\subsection{Coefficients for homology}  

All of our homology groups here, in whatever theory, are with
coefficients in $\Q$.  We expect things to work out similarly if we
used a different field as coefficients; it is less clear what would
happen if we tried to use $\Z$ as coefficients.

\subsection{Khovanov-Rozansky homology}

For the Khovanov-Rozansky homology, there are at least two separate
choices needed to fix a normalization.  The first is the normalization
of the HOMFLY polynomial, and the second is whether you want to view
the theory as homology or cohomology.  Most sources view it as
cohomology (e.g.~\cite{Khovanov, DBN}), but here we choose to view it
as a homology theory.  To make it a homology theory, we take the
standard cohomological chain complex and flip the homological grading
by $i \mapsto -i$, so that the differentials are now grading decreasing.
(One could also make it a homology theory by taking the dual complex
with dual differentials, but that's not what we do.)
  
For instance, to put a Poincar\'e polynomial $\khr_2(q,t)$ computed by
the Knot Atlas \cite{KnotAtlas} or KhoHo \cite{KhoHo} into our
conventions, one needs to substitute $q \to q^{-1}$ and $t \to t^{-1}$.  
(The first substitution is due to the differing conventions for
the Jones polynomial.)  Notice that this change has the same effect as
keeping the conventions fixed and replacing a knot by its mirror
image.

\subsection{Knot Floer homology} 

Our conventions for knot Floer homology $\HFK$ are opposite of the
usual ones in \cite{OShfk,Rasmussen}; in particular, our knot Floer
homology is the standard knot Floer homology of the mirror.  This has
the effect of simultaneously flipping both the homological and
Alexander gradings (see \emph{e.g.}~\cite[Eqn.~13]{OShfk}).  In
addition, we use different conventions for writing Poincar\'e
polynomials $\hfk$ than
\cite{KFC}.  For consistency with viewing
the Alexander polynomial $\Delta(K)$ as a specialization of the HOMFLY
polynomial, we view $\Delta(K)$ as a the polynomial in $q^2$ given by
$\Delta(K) = \pred(K)(a=1, q)$.  The variable $t$ in $\hfk$ gives the
homological grading.  In \cite{KFC}, $t$ is the variable for $\Delta(K)$
and $u$ is used for the homological grading; one can translate
information there into our conventions via the substitution $t \mapsto
q^{-2}, u \mapsto t^{-1}$.


\section{Families of differentials and relation to knot homologies}\label{sec-differentials}

As discussed in Section~\ref{sec-intro-redsuper}, we can only expect
uniform behavior for the $sl(N)$ homology for large $N$.  In this
section, we detail the additional structure that should encode the
$sl(N)$ homology for all $N$, and knot Floer homology as well.  We
start by presuming homology groups $\hk_{i,j,k}(K)$ categorifying the
reduced HOMFLY polynomial $\pred(K)(a,q)$.  The Poincar\'e polynomial
of this homology is the superpolynomial given by
\begin{equation}
 \Pred(K)(a,q,t)  = \sum  a^i q^j t^k \dim \hk_{i,j,k}(K).
\end{equation}
In addition, $\hk_*(K)$ should be equipped with a family of
differentials $\{d_N\}$ for ${N \in \Z}$, which will give the different
homologies.  The differentials should satisfy the following axioms:
\begin{list}{}
 {\setlength{\leftmargin}{1.3cm}
\setlength{\itemsep}{0.3cm}}
\item[\textbf{Grading:}] For \(N>0\),
\(d_N\) is triply graded of degree \((-2,2N,-1)\), {\it
  i.e.} 
\begin{equation*}
d_N: \hk_{i,j,k}(K) \to \hk_{i-2,j+2N,k-1}(K).
\end{equation*}
\(d_0\) is graded of degree \((-2,0,-3)\), and
 for \(N < 0\), \(d_N\) has degree \((-2,2N,
-1+2N)\).
\item[\textbf{Anticommutativity:}]  \(d_Nd_M = -d_Md_N\) for all \(N, M \in \Z\). In particular, \(d_N^2 = 0 \) for each \(N\in \Z\).
\item[\textbf{Symmetry:}] There is an involution \(\phi : \hk_{i,j,*} \to
  \hk_{i,-j,*} \) with the property that 
  \[
  \phi d_N = d_{-N} \phi \quad \mbox{for all \(N \in \Z\).}
  \]
\end{list} 
To build the connection to the other homology theories, first notice we
get a categorification of $P_N(K)$ by amalgamating groups to define
\begin{equation}
\hk^N_{p,k}(K) = \bigoplus_{iN+j= p} \hk_{i,j,k}(K).
\end{equation}
The Poincar\'e polynomial of these new groups is just $\Pred(K)(a = q^N,
q,t)$.  For $N > 0$, the first two axioms above imply that \(
(\hk^N_{l,k}(K),d_N) \) is a bigraded chain complex.  We can now state
our main conjecture:
\begin{conjecture}\label{Conj:Diffs}
  There is a homology theory $\hk_*$ categorifying the HOMFLY
  polynomial, coming with differentials $\{ d_N \}$ satisfying the
  three axioms.  For all $N>0$, the homology of $(\hk^N_*(K),d_N)$ is
  isomorphic to the $sl(N)$ Khovanov-Rozansky homology.  For $N = 0$,
  \((\hk^0_*(K),d_0)\) is isomorphic to the knot Floer homology. 
\end{conjecture}
For the last part of this conjecture, one must do
additional regrading of $\hk^0_*(K)$ to make it precise, see
Section~\ref{Subsec:HFK} below.  Let us denote the Poincar\'e polynomial
of the bigraded homology of $(\hk^N_*(K),d_N)$ by $\Pred_N(K)$; the
Khovanov-Rozansky part of the conjecture is thus summarized as
$\Pred_N(K) = \khr_N(K)$.

A few general comments are in order. First, for any given knot \(K\),
the superpolynomial has finite support, so the grading condition
forces \(d_N\) to vanish for \(N\) sufficiently large.  Thus the
earlier Conjecture~\ref{khredconj} is a special case of
Conjecture~\ref{Conj:Diffs}.

Second, we remark that the symmetry property generalizes the
well-known symmetry of the HOMFLY polynomial:
\begin{equation}
\pred(K)(a,q) = \pred(K)(a,q^{-1}). 
\end{equation}

 Finally, the homological grading of \(d_N\) for \(N<0\) may strike
 the reader as somewhat peculiar.  As we will explain in
 Section~\ref{Subsec:Delta},
  it is a natural consequence of the symmetry \(\phi\).


\subsection{Examples} 

To illustrate the properties above, we consider three examples,
starting with the easy case of the unknot. 

\begin{example:unknot}
  For the unknot $U$, all the $sl(N)$ homology is known and $\khr_N(U) =
  1$ for all $N > 0$.  Thus the superpolynomial is clearly given by \(\CP (U) = 1\),
  where all the differentials \(d_N\) are identically zero.
\end{example:unknot}

\begin{example:T23}
 The HOMFLY polynomial of the negative
trefoil knot $T_{2,3}$ is given by \(P (T_{2,3}) = a^2q^{-2}+a^2q^2 - a^4\). The
corresponding superpolynomial also has three terms:
\begin{equation}
\CP (T_{2,3}) = a^2q^{-2}t^0+a^2q^2t^2 + a^4q^0t^3
\end{equation}
To illustrate the differentials, it is convenient to represent 
\(\hk(K)\) by a \emph{dot diagram} as shown in Figure~\ref{Fig:Tref}. 
\begin{figure}

\begin{center}
\begin{minipage}[b]{15cm}
\parbox[c]{7.5cm}{
\begin{overpic}[tics=10]{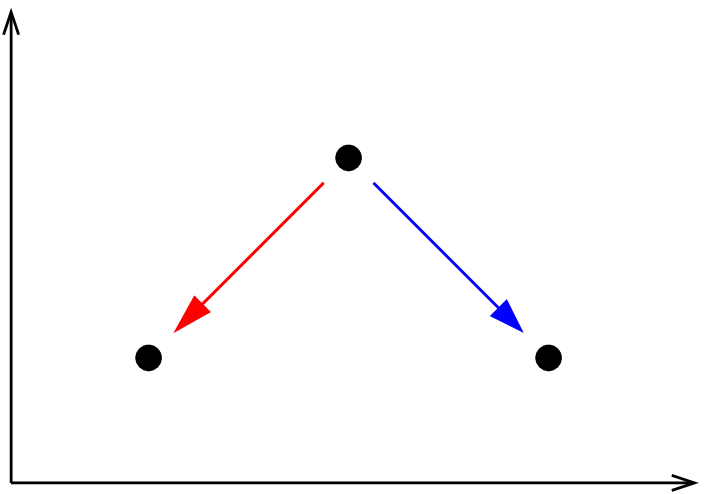}
  \put(13, 9){$a^2 q^{-2} t^0$}
  \put(72, 9){$a^2 q^2 t^2$}
  \put(43, 53){$a^4 q^0 t^3$}
  \put(26, 37){$d_{-1}$}
  \put(64, 37){$d_{1}$}
  \put(-4, 65){$a$}
  \put(100, 0){$q$}

\end{overpic}}
\hfill
\parbox[c]{5cm}{%

\

\vspace{0.5cm}

  \includegraphics{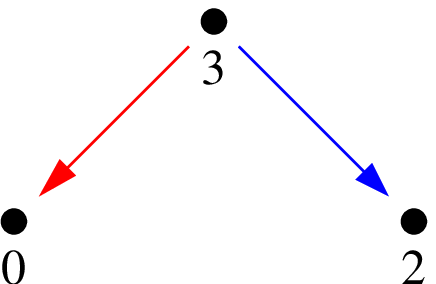}

}
\end{minipage}
\end{center}
\vspace{0.5cm}

\caption{\label{Fig:Tref} 
  Non-zero differentials for the trefoil knot.  On the left is a fully
  labeled diagram, and on the right is the more condensed form
  that we will use from now on.  The minimum $a$-grading is $2$.}
 
\end{figure}
We draw one dot for each term in the superpolynomial, so that the
total number of dots is equal to the dimension of \(\hk (K)\). The
dots' position on horizontal axis records the power of \(q\), and on
the vertical, the power of \(a\). The left-hand side of
Figure~\ref{Fig:Tref} shows such a diagram for the trefoil, with each
dot labeled by its corresponding monomial. 

Since the relative \(a\) and \(q\) gradings are determined by the
position of the dots, we omit them from the diagram and just label
each dot by its \(t\)-grading.  To fix the absolute \(a\)-grading, we
record the \(a\)-grading of the bottom row.  Determining the absolute
$q$-grading from such a picture is easy, since the line \(q=0\)
corresponds to the the vertical axis of symmetry.
The nonzero components of \(d_i\) are shown by arrows of slope
\(-1/i\).  As indicated by the figure, the trefoil has two nontrivial
differentials: \(d_1\) and \(d_{-1}\).

Now let's substitute \(a=q^N\) and take homology with respect to
\(d_N\).  For \(N>1\), there are no differentials, and so we just get
\(\Pred_N(T_{2,3}) = q^{2N-2}t^0 + q^{2N+2}t^2 + q^{2N}t^3 \).  For
\(N=1\), the differential \(d_1\) kills the two right-hand generators,
and we are left with \(\Pred_1(T_{2,3}) = 1\). In this case, it is
possible to check directly that $\Pred_N = \khr_N$ for all \(N>0\).
Note that $\khr_1$ of any knot is always $1 = q^0t^0$, which is why $d_{\pm 1}$
must be non-zero even in such a simple example as this.
\end{example:T23}

\begin{example:T34}
  A more complicated example is provided by the negative \((3,4)\)
  torus knot, which is the mirror of the knot $8_{19}$.  In this case,
  both the HOMFLY polynomial and the superpolynomial have \(11\)
  nontrivial terms:
  \begin{equation}
  \begin{split}
  \pred(T_{3,4}) &= a^{10} - a^8(q^{-4} + q^{-2}  +1  + q^2 + q^4) +  a^6(q^{-6} + q^{-2} + 1  + q^2 + q^6) \\ 
  \Pred(T_{3,4}) &= a^{10}t^8 + a^8(q^{-4}t^3 + q^{-2}t^5  +t^5  + q^2 t^7 + q^4 t^7) +  a^6(q^{-6} t^0 + q^{-2}t^2 + t^4  + q^2 t^4 + q^6t^6)
 \end{split}
 \end{equation}
 The superpolynomial is illustrated by the dot
  diagram in Figure~\ref{Fig:T34}.
\begin{figure}
\centering\includegraphics{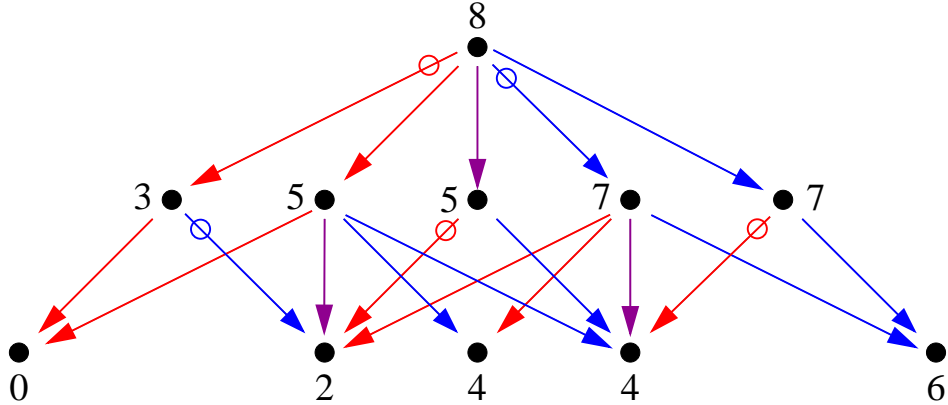}
\caption{\label{Fig:T34} Differentials for \(T_{3,4}\).The
  bottom row of dots has \(a\)-grading \(6\).  The leftmost dot on
  that row has $q$-grading $-6$, which you can determine by noting
  that the vertical axis of symmetry corresponds to the line \(q=0\). }
\end{figure}

Here
there are five nontrivial differentials: \( d_{-2},d_{-1},d_0,d_1\),
and \(d_2\).  To understand the differentials completely, think of the
dots as representing specific basis vectors for $\hk_{i,j,k}$; then an
arrow means the corresponding $d_N$ takes the basis element at its
tail to $\pm$ the basis element at its tip.  In this case, the sign can
be inferred from the diagram; those that switch the sign have small
circle at their tails.  (To avoid clutter, in future we will leave it
to reader to choose appropriate signs for the differentials.)  It is
now easy to check that all the $d_N$ anticommute.  The symmetry
involution $\phi$ corresponds to flipping the diagram about its
vertical axis of symmetry.  For the $\hk_*$ off the line itself, $\phi$
permutes our preferred basis vectors; on $\hk_{10, 0, 8}$ and
$\hk_{8,0,5}$ it acts by $-\mathrm{Id}$, but is the identity on
$\hk_{6,0,4}$.  You can now easily check the symmetry axiom.  

Substituting \(a=q^2\) and taking homology with respect to \(d_2\) kills
six generators, leaving 
\begin{equation}
\Pred_2(T_{3,4}) = q^6t^0 + q^{10}t^2 +q^{12}t^3+q^{12}t^4+q^{16}t^5
\end{equation}
which is the ordinary (\(N=2\)) Khovanov homology of \(T_{3,4}\).  As
before, $\Pred_1(T_{3,4}) = 1$ --- only the bottom leftmost term survives.
\end{example:T34}


\subsection{Relation to knot Floer homology}
\label{Subsec:HFK}

In order to recover the knot Floer homology, we must introduce a new
homological grading on \(\hk(K)\), which is given by \(t'(x) = t(x) -
a(x)\). In other words, the Poincar\'e  polynomial of \(\hk\) with
respect to the new grading is
\begin{equation}
\Pred'(a,q,t) = \Pred(a = at^{-1},q,t).
\end{equation}
The differential \(d_0\) lowers the new grading \(t'\) by \(1\).  Now
forget the $a$-grading (\emph{i.e.}~substitute \(a=1\)), and take the
homology with respect to \(d_0\).  We denote the Poincar\'e polynomial
of this homology by $\Pred_0(K)(q, t)$, and this homology categorifies
the Alexander polynomial $\Delta(K)(q^2) = \pred(K)(a = 1, q)$.  A precise
statement of the last part of Conjecture~\ref{Conj:Diffs} is that
$\Pred_0(K) = \hfk(K)$, where $\hfk$ is the Poincar\'e polynomial of
knot Floer homology defined in \eqref{hfkdef}.

\begin{figure}
\centering\includegraphics{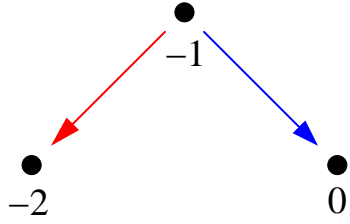}
\caption{\label{Fig:Regrade} Trefoil with new homological gradings.}
\end{figure}

As a first example of this process, consider the trefoil
knot. Figure~\ref{Fig:Regrade} shows the generators for
\(\hk(T_{2,3})\) with respect to the new homological grading
\(t'\). The differential \(d_0\) is trivial, so we find 
\begin{equation}
\Pred_0(T_{2,3}) = \Pred(T_{2,3})(a=t^{-1}, q, t) = q^{-2}t^{-2} + q^0t^{-1}+q^2t^0,
\end{equation} 
which is indeed equal to \(\hfk(T_{2,3})\).  

Next we consider \(T_{3,4}\), for which \(d_0\) kills \(6\) of the
\(11\) generators. We leave it to the reader to check that after
regrading and taking homology with respect to \(d_0\), we are left
with
\begin{equation}
\Pred_0(K) = q^{-6}t^{-6}+q^{-4}t^{-5}+q^0t^{-2}+q^4t^{-1}+q^6t^0,
\end{equation}
which agrees with $\hfk(T_{3,4})$.  

\subsection{The $\delta$-grading and symmetry}
\label{Subsec:Delta}

It is natural to consider a fourth grading on $\hk(K)$ which is
obtained as a linear combination of the \(a,q,\) and \(t\) gradings.
It is defined by
\begin{equation}
\delta(x) = t(x) - a(x) - q(x)/2.
\end{equation}
When we specialize to  \(\HFK\) or \(\hkr_2\), the
\(\delta\)-grading reduces to the \(\delta\)-gradings on these two
theories defined in \cite{Rasmussen}. 
Indeed, if \(q_2\) is the \(q\)-grading on \(\hkr_2\) defined by setting
\(a=q^2\), then 
\begin{equation}
 t(x) - a(x) - q(x)/2 = t(x) - (2a(x)+q(x))/2 = t(x) - q_2(x)/2,
\end{equation}
where \(q_2\) denotes the \(q\)-grading on \(\hkr_2\) and 
the rightmost expression is the \emph{definition} of the
$\delta$-grading on $\hkr_2$.  Similarly, if \(t'\) is the homological
grading on \(\HFK\), defined by setting \(a = 1/t\), then
\begin{equation}
t(x) - a(x) - q(x)/2 = t'(x) - q(x)/2
\end{equation}
where the right-hand side is the definition of the \(\delta\)-grading
on \(\HFK\). 

We can use the \(\delta\)-grading to justify the somewhat peculiar
behavior of \(d_{i}\)  for \(i<0\) with respect to the homological
grading. In analogy with knot Floer homology, where the \(\delta\)-grading
is preserved by the conjugation symmetry, we expect that the
\(\delta\)-grading will be preserved by the symmetry \(\phi\) of
Conjecture~\ref{Conj:Diffs}. For \(i>0\), the differential \(d_i\) lowers the
\(\delta\)-grading by \(1-i\). Since \( \phi\) exchanges \(d_i\) and
\(d_{-i}\), the differential \(d_{-i}\) should lower the \(\delta\)
grading by \(1-i\) as well. It is then easy to see that \(d_{-i}\)
lowers the homological grading by \(-1-2i\).


\subsection{Canceling differentials}
\label{Subsec:CDiffs}
Let \((C,d)\) be a chain complex. We say that \(d\) is a {\it
  canceling differential} on \(C\) if the homology of \(C\) with
  respect to \(d\) is one-dimensional. The presence of a canceling
  differential is an important feature of all the reduced knot
  homologies. For \(\HFK\), this was known from the start ---
  essentially, it's the fact that \(\hfhat(\Sp^3) \cong \Z\). 
For the $sl(2)$ Khovanov homology, it follows from work of Turner
  \cite{Turner}, which itself builds on work of Lee \cite{ESL2} and
  Bar-Natan \cite{DBN04}. Finally, the existence
  of such a differential for \(\hkr_N\) can be derived by combining Turner's results
  with the work of Gornik \cite{Gornik} in the unreduced case. 

Conjecture~\ref{Conj:Diffs} provides a unified explanation for the
presence of these canceling differentials. Indeed, for any knot \(K\),
 \(\CP_1(K)= 1\), which implies that 
\(d_1\) should be  a canceling differential on \(\hk(K)\). We expect that the
known differentials on the various specializations of \(\hk\) are all
induced by the action of \(d_1\). 

To state this more precisely, let us suppose that
Conjecture~\ref{Conj:Diffs} is true. Since \(d_1\) anticommutes with
\(d_N\), the pair \((\hk(K),d_1+d_N)\) is also a chain
complex. Consider the grading on \(\hk(K)\) obtained by setting
\(a=q^N\). This grading is preserved by \(d_N\), but is strictly
lowered by \(d_1\). In other words, it makes \((\hk(K),d_1+d_N)\) into
a filtered complex whose associated graded complex is
\((\hk(K),d_N)\). Since we are using rational coefficients, we can
 reduce this complex to a chain homotopy equivalent complex of the
form \((H_*(\hk(K),d_N),d_1')\). (See Lemma 4.5 of
\cite{Rasmussen} for a proof.)

\begin{prop}
If we assume that Conjecture~\ref{Conj:Diffs} holds, then 
\(d_{1}'\) is a canceling differential on \( H_*(\hk(K),d_N)\) whenever
$N \neq 1$. 
\end{prop}

\begin{proof}
We again consider the complex \((\hk(K),d_1+d_N)\), but with a
different grading --- namely, the one defined by setting \(a=q\). It is
easy to see that \(d_1\) preserves the new grading, while \(d_N\)
strictly raises it, so this grading also makes \( (\hk(K),d_1+d_N)\)
into a filtered complex. Reducing as before, we obtain a chain
homotopy equivalent complex \((H_*(\hk(K),d_1),d_N')\). Assuming that
the conjecture is true, \(H_*(\hk(K),d_1) \cong \hkr_1 (K)\)
is one-dimensional, so
\begin{align}
H_*(H_*(\hk(K),d_N),d_1') & \cong H_*(\hk(K),d_1+d_N) \\
& \cong  H_*(H_*(\hk(K),d_1),d_N') \\ 
& \cong H_*(\hk(K),d_1)
\end{align}
is one-dimensional as well. 
\end{proof}

An interesting consequence of Conjecture~\ref{Conj:Diffs} is that it
predicts the existence of a {\it second} canceling differential on
\(\hkr_N\). Indeed, the symmetry property implies that
 \(d_{-1}\) is also a canceling differential on \(\hk\), and the same
 argument used for \(d_1\) implies that it should descend to a differential
 on any specialization of \(\hk\). 
 
 In the case of \(\HFK\), it is well-known that two such differentials
 exist, and that they are exchanged by the conjugation symmetry (see,
 {\it e.g.}~\cite[Prop.~4.2]{Rasmussen}).  To illustrate this fact, we
 consider the knot Floer homology of the trefoil.  There,
 \(\HFK(T_{2,3})\) has three generators, corresponding to monomials
 \(q^{-2}t^{-2}\), \(q^0t^{-1}\), and \(q^2t^0\) in the Poincar{\'e}
 polynomial. Looking at Figure~\ref{Fig:Regrade}, we see that
 differential induced by \(d_{-1}\) takes the second generator to the
 first , while the differential induced by the \(d_1\) takes the
 second generator to the third. This is indeed the differential
 structure on \(\HFK(T_{2,3})\).
 
 In general, the differential induced by \(d_{-1}\) should correspond
 to the usual differential on \(\HFK\) (that is, the one that lowers
 the Alexander grading), while the differential induced by \(d_{1}\)
 corresponds to its conjugate symmetric partner. As a check, let us
 consider how the two induced differentials behave with respect to the
 homological grading \(t'\). Since both \(d_0\) and \(d_{-1}\) lower
 the homological grading by \(1\), the induced map \(d_{-1*}\) will
 lower \(t'\) by \(1\) as well. This is in accordance with the
 behavior of the usual differential on \(\HFK\). In contrast, \(d_1\)
 raises \(t'\) by \(1\), so the behavior of \(d_{1*}\) with respect
 to \(t'\) is somewhat more complicated. In fact, it is not hard to
 see that if some component of \(d_{1*}\) raises the \(q\)-grading by
 \(2k\), it will raise \(t'\) by \(2k-1\). This is precisely the
 behavior exhibited by the ``conjugate'' differentials in knot Floer
 homology.

In contrast, the differential \(d_N\) which  gets us from
\(\hk(K)\) to \(\hkr_N(K)\) lowers the usual homological grading on
\(\hk(K)\) by
\(1\), does \(d_1\). Thus the 
 differential induced by \(d_1\) on \(\hkr_N(K)\) will respect the
 homological grading on that group.
 We expect that \(d_{1*}\) corresponds to the
 differential of Lee, Turner and Gornik. As an example
 consider the \(sl(2)\) homology of the trefoil. 
Here, we have \(\CP_2(T_{2,3})
= q^2t^0+q^6t^2 + q^8t^3\), and the differential induced by \(d_1\)
takes the third term to the second. This agrees with the standard
canceling differential on the reduced Khovanov homology.

As far as we are aware, the presence of a second canceling
differential on the Khovanov homology has not been considered before.
Although we do not know how to construct such a differential directly,
in section~\ref{Subsec:d-1} we describe some evidence which supports
the idea that \(\hkr_2\) admits an additional canceling differential
induced by \(d_{-1}\). 


\subsection{Analog of $s$ and  $\tau$}

Given a canceling differential on a filtered chain complex, one can
define a simple invariant by considering the filtration grading of the
(unique) generator on homology. Applying this fact to knot Floer
homology, Ozsv{\'a}th and Szab{\'o} \cite{OStau} 
defined a knot invariant \(\tau(K)\), which carries
information about the four-ball genus of \(K\). Subsequently, an analogous
invariant \(s\) was defined using the Khovanov homology  \cite{khg}. 

On the triply graded homology theory \(\hk(K)\), the canceling
differential \(d_1\) can be used to define a similar
invariant. Since there are two polynomial gradings on \(\hk(K)\), it
initially looks  like we will get two invariants. In reality,
however, the generator of the homology with respect to \(d_1\) 
always lies on the line where \(q(x) = -a(x)\). This is because when we
specialize to the $sl(1)$ theory by substituting \(a=q\), the generator
corresponds to the unique term in \(\CP_1(K) = 1\). After taking
homology with respect to \(d_1\),  the surviving term will have
the form \(a^Sq^{-S}t^0\). The number \(S\) will be an invariant of
\(K\) analogous to \(s\) and \(\tau\). 

For example, if \(K\) is the \((3,4)\) torus knot, a glance at
Figure~\ref{Fig:T34} shows that \(S(K)=6\). This example illustrates
an interesting feature of \(S\): namely, that it is in some sense  
easier to compute than either \(s\) or \(\tau\). Indeed, to compute
\(S\), we need only consider those generators of \(\hk(K)\) which lie
along the line \(a(x) = -q(x)\). In many cases (like the one above)
the number of generators we need to consider is quite small. 

In analogy with the known properties of \(S\) and \(\tau\), we expect
that \(S\) will be a lower bound for the four-ball genus of \(K\) (see
Section~\ref{Subsec:four-ball}). It is not clear, however, whether it
contains any new information, since in all the examples we have
considered, it appears that \(S(K)=s(K)=2\tau(K)\). We hope that further
consideration of the construction of \(S\) will shed new light on the
relationship between \(s\) and \(\tau\), either by proving that all
three quantities are equal, or by suggesting where to look for a
counterexample.


\subsection{Motivation for the conjecture} 

We conclude this section by briefly sketching the background to
Conjecture~\ref{Conj:Diffs}, and indicating how strongly we believe
its various parts. Our main reason for expecting the presence of the
differentials \(d_N\) for \(N>0\) comes from Gornik's work on the
\(sl(M)\) homology. In \cite{Gornik}, Gornik describes a deformation
of Khovanov and Rozansky's construction which gives rise to a
canceling differential on \(\hkr_M\). In fact, this construction may
be easily modified to obtain a whole family of deformations, one for
each monic polynomial of degree \(M\). It follows that  any monic
polynomial of degree \(M\) gives rise to a differential on \(\hkr_M\).
If we let \(d_N^{(M)}\) be the differential corresponding to the
polynomial \(X^M-X^N\), we expect the differential \(d_N\) of the
conjecture can be obtained as the limit of \(d_N^{(M)}\) as \(M\to
\infty\).   In analogy with Gornik's work,
we expect that taking the homology of \(\hkr_M(K)\) 
with respect to this differential \(d_N^{(M)}\) will 
 give the group  \(\hkr_N(K)\), thus matching the behavior predicted
by Conjecture~\ref{Conj:Diffs}. (Indeed, this observation was the
genesis of the conjecture.) 
For \(N>0\), the behavior expressed by the grading axiom was chosen to
agree with the known behavior of \(d_N^{(M)}\). 
Finally, the fact that \(d_{N_1}\) and \(d_{N_2}\) \((N_1,N_2>0)\) 
anticommute  should follow from the 
  linearity of the space of deformations. More precisely, if we let
\(d_{N_1,N_2}^{(M)}\) be the differential corresponding to the polynomial
\(X^M-X^{N_1}-X^{N_2}\), then 
\(d_{N_1,N_2}^{(M)} = d_{N_1}^{(M)}+d_{N_2}^{(M)}\), so
the fact that \((d_{N_1,N_2}^{(M)})^2=0\) implies that \(d_{N_1}^{(M)}\)
    and \(d_{N_2}^{(M)}\) anticommute.

The rest of the conjecture is more speculative. Our original reason
for expecting the presence of the differentials \(d_N\) for \(N \leq
0\) was based on analogy with the knot Floer homology. We believe that
the strong internal consistency of the theory, as seen in
the examples of Section~\ref{Sec:Examples}, together with
the apparently correct predictions it makes (such as the computations
of the stable $sl(2)$ Khovanov homology of the torus knots in
Section~\ref{Subsec:KHR2})  indicate that there must be 
{\it something} meaningful going on. It is possible, however,
that  we have erred in stating the exact details. Below, we outline 
some potential weak points of Conjecture~\ref{Conj:Diffs}.

\begin{itemize}
\item We are not currently aware of any
construction which might give rise to the \(d_N\)'s for \(N \leq 0\).
Our reasons for expecting their existence are based on analogy with
the case \(N>1\), which suggests that there should be a differential
\(d_0\) giving rise to knot Floer homology, and with knot
Floer homology itself, whose symmetries suggest the presence of
\(d_N\) for \(N <0\).

\item The statement in the conjecture about the gradings of
differentials is somewhat stronger than would be expected from
Gornik's work. 
{\it A priori}, the differentials coming from Gornik's theory should
shift the \((a,q)\) bigrading by some multiple of \((-2,2N)\).
 The requirement that this multiple is always one is imposed to ensure
 that \(d_N\) shifts both \(t\) and \(t'\) by a constant amount. 
(Some further
support for this idea is provided by the fact that 
there are a number of ten-crossing knots which at first glance
look as if \(d_{1}\) might lower the \((a,q)\) bigrading by
\((-4,4)\). In all these examples, however, further examination
suggests that this is not the case.)

\item 
Finally, there is some chance that taking homology with respect to \(d_0\)
does not give the knot Floer homology, but some other
categorification of the Alexander polynomial which  happens to
look a lot like it. An interesting test case
for this possibility is provided by the presence of 
mutant knots with different genera. For example, there are several
mutant pairs of \(11\)-crossing knots, one of which has genus one
bigger than the other. These knots have the same HOMFLY polynomial and
\(\khr_2\), but their knot Floer homologies must differ. It is an
interesting question to determine whether these knots have the same
superpolynomial and (if they do) the same differentials. 
\end{itemize}


\section{Geometric interpretation}
\label{sec-geometry}

In this section, we explain in more detail the geometric
interpretation of the triply graded knot homology in the language of
open Gromov-Witten theory.  As discussed in
Section~\ref{subsec-intro-super}, this relation was part of the
original motivation for the triply graded theory, and we hope it can be
useful for developing both sides of the correspondence.
In this section, we mainly consider the unreduced homology
which has a more direct relation to the geometry of holomorphic curves.

The geometric setup consists of the following data:
a non-compact Calabi-Yau 3-fold $X$ and a Lagrangian
submanifold $\CL \subset X$. Therefore, for every
knot $K \subset \Sp^3$, we need to define $X$ and $\CL$.
The Calabi-Yau space $X$ is independent of the knot;
it is defined as the total space of
the $\CO (-1) \oplus \CO (-1)$ bundle over $\cp^1$:
\begin{equation}
\label{oobundlex}
\CO (-1) \oplus \CO (-1) \to \cp^1
\end{equation}
On the other hand, the information about the knot $K$
is encoded in the topology of the Lagrangian submanifold,
which we denote $\CL_K$ to emphasize that it is determined
by the knot:
\begin{equation}
K \leadsto \CL_K
\end{equation}
A systematic construction of the Lagrangian submanifold $\CL_K$ from
a braid diagram of $K$ was proposed by Taubes \cite{Taubes}.
It involves two steps. First, one constructs a two-dimensional
non-compact Lagrangian submanifold $\CL_K^{(2)} \subset \C^2$,
which has the property that its intersection with a large
radius 3-sphere, $\Sp^3 \subset \C^2$, is isotopic to the knot $K$.
Then, we identify $\C^2 \otimes \CO (-1)$ with a fiber of $X$
and define $\CL_K$ to be a particular subbundle
$\CL_K^{(2)} \to \Sp^1$ of the bundle \eqref{oobundlex}
restricted to the equator $\Sp^1 \subset \cp^1$.
The construction is such that $\CL_K$ is Lagrangian
with respect to the standard K\"ahler form on $X$.
Moreover, for every knot $K$, the resulting 3-manifold $\CL_K$
has the first Betti number $b_1 (\CL_K)=1$.

Given a Calabi-Yau space $X$ and a Lagrangian submanifold
$\CL_K \subset X$, it is natural to study holomorphic Riemann
surfaces in $X$ with Lagrangian boundary conditions on $\CL_K$:
\begin{equation}
\label{embeddings}
(\Sigma,\partial \Sigma) \hookrightarrow (X,\CL_K)
\end{equation}
Specifically, we consider embedded surfaces $\Sigma$
which satisfy the following conditions:
\begin{enumerate}

\item $\Sigma$ is a holomorphic Riemann surface with a fixed genus $g$
and one boundary component, $\partial \Sigma \cong \Sp^1$,

\item $[\Sigma] = Q$ with $Q$ a fixed class in $H_2 (X, \CL_K ;\Z) \cong \Z$,

\item $[\partial \Sigma] = \gamma$, where $\gamma$ 
 generates  the free part of the homology group
$H_1 (\CL_K, \Z) \cong \Z \gamma$ (mod torsion).

\end{enumerate}
Now we are ready to define the moduli spaces that appear
in the geometric interpretation of the triply graded theory,
{\it cf.} \eqref{hmgq}.
Let $\Sigma$ be an embedded Riemann surface
which satisfies the conditions (1) -- (3),
and let $A \in \Omega^1 (\Sigma)$
be a flat $U(1)$ gauge connection on $\Sigma$,
\begin{equation}
F_A = 0.
\end{equation}
We define $\CM_{g,Q} (X,\CL_K)$ to be moduli ``space''
of the embedded Riemann surfaces $\Sigma$ with a gauge connection $A$,
modulo the gauge equivalence, $A \to A + df$ where $f \in \Omega^0 (\Sigma)$.
Assuming that the dependence on $X$ and $\CL_K$ is clear from
the context, we often refer to this moduli space simply as $\CM_{g,Q}$.
The cohomology groups $H^k (\CM_{g,Q})$ are labeled by three integers:
the degree $k$, the genus $g$, and the relative homology class
$Q \in H_2 (X, \CL_K ;\Z) \cong \Z$. These are the three gradings
of our triply graded theory.

\begin{remark}
Since in general $\CM_{g,Q}$ may be singular and non-compact,
one needs to be careful about the definition of $H^k (\CM_{g,Q})$.
This problem is familiar in the closely related context
of Gromov-Witten theory, where instead of embedded Riemann
surfaces with a flat connection one ``counts'' stable
holomorphic maps (possibly with boundary).
In Gromov-Witten theory, there is a way to define cohomology
classes and intersection theory on the moduli spaces of stable maps
(see \cite{KatzLiu,LiSong,GraberZ} for some recent work on
the mathematical formulation and calculation of the open
Gromov-Witten invariants).
Similarly, the physical interpretation of the $sl(N)$ knot
homology \cite{GSV} suggests that, at least in the present case,
there should exist a suitable definition of $\CM_{g,Q}$,
such that the cohomology groups $H^k (\CM_{g,Q})$ can be
identified with the triply graded knot homology groups.
\end{remark}

\begin{example:unknot}
  In this case, the only non-trivial holomorphic curves are
  holomorphic disks wrapped on the northern and the southern
  hemispheres of the $\cp^1 \subset X$.  Their moduli spaces are isolated
  points, $\CM_{g,Q} \cong \mathrm{pt}$ for $g=0$ and $Q = \pm 1$, which correspond
  to the two terms, $a$ and $a^{-1}$, in the unreduced superpolynomial
  for the unknot $\Pun(a, q, t) = a - a^{-1}$.
\end{example:unknot}

\subsection{Genus expansion and symmetry}

Now, let us look more closely at the structure of the moduli
space $\CM_{g,Q}$, assuming that it is well defined.
Let $\Sigma$ be a non-degenerate Riemann surface of genus $g$.
The moduli space of gauge equivalence classes of flat
$U(1)$ connections $A \in \Omega^1 (\Sigma)$ is
isomorphic to a $2g$-dimensional torus,
\begin{equation}
\mathrm{Hom} (\pi_1 (\Sigma);U(1))/ U(1) \cong T^{2g}.
\end{equation}
Therefore, $\CM_{g,Q}$ has the structure of a fibration
\begin{equation}
\label{mspacefibr}
\begin{array}{ccc}
T^{2g} & \to & \CM_{g,Q} \\
&& \downarrow \\
&& \CM^{\mathrm{geom}}_{g,Q}
\end{array}
\end{equation}
where $\CM^{\mathrm{geom}}_{g,Q}$ is the moduli space
of embedded Riemann surfaces \eqref{embeddings}
which satisfy the conditions (1) -- (3).
In many cases, the fibration structure \eqref{mspacefibr}
can be recognized directly in the structure of the superpolynomial
written in terms of the variables $a$, $t$, and $y$,
where $y = (qt^{1/2} + q^{-1} t^{-1/2})^2$.
In particular, the contribution of an isolated Riemann surface
with genus $g$ and relative homology class $Q$ looks like \cite{GSV}:
\begin{equation}
a^Q t^r (qt^{1/2} + q^{-1} t^{-1/2})^{2g}
\end{equation}
where the last factor is the familiar Poincare polynomial of $T^{2g}$.
In general, the superpolynomial $\Pun (K)$ should have the following structure
\begin{equation}
\label{gexpansion}
\Pun (K) = \sum_{g,Q,i} \hat D_{Q,g,i} a^Q t^i (qt^{1/2} + q^{-1}t^{-1/2})^{2g}
\end{equation}
where $\hat D_{Q,g,i} \in \Z$ encode the geometry
of the fibration \eqref{mspacefibr}.
We refer to the expansion \eqref{gexpansion} as the genus expansion.
It is natural to expect a similar structure also in the case of
the reduced superpolynomial, $\Pred (K)$.
Notice, in the reduced case, the expansion of the form \eqref{gexpansion}
is equivalent to the existence of the symmetry
\begin{equation}
\label{symmphi}
\phi: \hk_{i,j,*} (K) \to \hk_{i,-j,*} (K)
\end{equation}
that we discussed earlier in Section \ref{sec-differentials}.
In the geometric interpretation, this symmetry follows from
the fibration structure \eqref{mspacefibr}.

For the genus expansion of the reduced superpolynomial,
let us also define the {\it holomorphic genus}, $g_h(K)$,
to be the maximum value of $g$ which occurs in the sum \eqref{gexpansion}.
It has a clear geometric meaning as the maximum genus
of the holomorphic Riemann surface \eqref{embeddings}
which satisfies the conditions (1) -- (3).
With this definition, $2g_h(K)$ is equal to the maximum
power of $q$ that appears in the reduced superpolynomial.
The conjectured relation with knot Floer homology suggests the following bound
\begin{equation}
\label{ghbound}
g_3 (K) \leq g_h (K)
\end{equation}
where $g_3(K)$ is the Seifert genus of $K$.

\subsection{Relation to Gromov-Witten invariants}

Let us conclude this section by noting that taking
the Euler characteristic in the triply graded knot
homology $\hk_* (K)$ translates into taking the Euler
characteristic in $H^* (\CM_{g,Q})$.
On the other hand, the invariants $\chi (\CM_{g,Q})$,
which in the physics literature are called ``integer BPS invariants'',
contain the information about all-genus open Gromov-Witten
invariants of $(X,\CL_K)$ \cite{OV,LMV}.
The relation between the open Gromov-Witten invariants
and the integer BPS invariants is very non-trivial.
For example, the genus-counting parameter $u$
in the open Gromov-Witten theory is related to
the variable $q$ that we use via the following change
of variables (also familiar in the context
of the closed Gromov-Witten theory \cite{MNOP}):
\begin{equation}
\label{qviau}
q = e^{iu}
\end{equation}
Via this relation, all the information about the relative
Gromov-Witten theory of $(X,\CL_K)$ can be compactly recorded
in a finite set of non-zero integer BPS invariants.
One can use this relationship both ways. In particular,
one can find the Euler characteristic $\chi (\CM_{g,Q})$ by
computing the open Gromov-Witten invariants, say via
the localization technique \cite{GraberZ,KatzLiu,LiSong}.
It would be interesting to extend the existing techniques to compute
the dimensions of the individual cohomology groups $H^k (\CM_{g,Q})$.


\section{Examples and patterns}
\label{Sec:Examples}

We now describe the superpolynomials associated to some specific knots
with \(10\) or fewer crossings.  Although we lack a definition for the
triply graded theory and are unable to compute the $sl(N)$ homology in
general, we can still make intelligent guesses at the form of the
superpolynomial, based on Conjecture~\ref{Conj:Diffs} and the known
values of \(\HFK\) and \(\hkr_2\).  These example illustrate the
internal consistency of the structure proposed in
Conjecture~\ref{Conj:Diffs}.  Once we have looked at these examples,
we explore some patterns observed there in more detail in
Sections~\ref{Subsec:ranks}--\ref{Subsec:d-1}.


\subsection{Thin knots}
\label{SubSec:Thin}

In both knot Floer homology and $sl(2)$ Khovanov homology, the
smallest knots exhibit the following simple behavior: if we plot the
homological grading versus the polynomial grading, all the generators
line up along a single line.  Moreover, this line always has the same
slope, which corresponds to the appropriate $\delta$-grading being
constant (see Section~\ref{Subsec:Delta} for definitions).  Such
knots are called \emph{thin} (with respect to either \(\HFK\) or
\(\hkr_2\)).  In the triply graded case, we can define thinness analogously:

\begin{definition} A knot \(K\) is \emph{\(\hk\)--thin} if all
generators of \(\hk(K)\) have the same \(\delta\)-grading. 
\end{definition}

For an $\hk$--thin knot, the $t$-grading of a term of $\Pred(K)$ is
determined by the $a$-  and $q$- gradings.  Thus, there  can be no
cancellation when we specialize $\Pred(K)$ to $\pred(K)$, and so
$\Pred(K)$ is completely determined by its HOMFLY polynomial and the
common \(\delta\)-grading of its generators.  Noting that the common
$\delta$-grading is equal to \(-S(K)/2\), the precise relationship
between $\Pred(K)$ and $\pred(K)$ is concisely expressed by:
\begin{equation}
\Pred_K(a,q,t) = (-t)^{-S(K)/2} P_K(at, iqt^{1/2}).
\end{equation}

If \(K\) is thin, the dimension of \(\hk(K)\) is equal to the determinant
 of \(K\). Moreover, all differentials other than \(d_1\) and
 \(d_{-1}\) automatically vanish, since these differentials lower the
 \(\delta\)-grading. Finally, the fact that \(d_1\) and \(d_{-1}\) anticommute
and each have one-dimensional homology implies that \(\hk(K)\) can be
 decomposed as the direct sum of a number of ``squares'' with
 Poincar{\'e} polynomial
 \(a^{i}q^jt^k(1+a^{-2}q^2t^{-1})(1+a^{-2}q^{-2}t^{-3})\) and a single
 ``sawtooth'' summand isomorphic to \(\hk(T_{2,k})\) for some value of
 \(k\). It follows that 
\begin{equation}
\Pred_K(a,q,t) = \Pred_{T_{2,k}}(a,q,t) + 
(1+a^{-2}q^2t^{-1})(1+a^{-2}q^{-2}t^{-3})Q(a,q,t)
\end{equation}
where \(Q\) is a polynomial with positive coefficients. We thus obtain
a restriction on the HOMFLY polynomial of a thin knot: if \(T_{2,k}\)
is a torus knot whose signature is equal to \(S(K)\), the polynomial
\begin{equation}
\frac{P(K) - P(T_{2,k})}{(1-a^{-2}q^2)(1-a^{-2}q^{-2})}
\end{equation}
must be alternating.

 As with \(\HFK\) and \(\hkr_2\), we expect some classes of simple
 knots are \(\hk\)--thin.  In particular,

\begin{conjecture}\label{Conj:two-bridge}
If \(K\) is a two-bridge knot, then \(K\) is \(\hk\)--thin, and
\(S(K)=\sigma(K)\).
\end{conjecture}

As two-bridge knots are alternating and hence thin for \(\HFK\) and
\(\hkr_2\) \cite{OSAlt, ESL1}, it is easy to check that
Conjecture~\ref{Conj:two-bridge} holds for $N = 0,1,2$.  Thus, to
prove it one needs to show
\begin{equation}\label{eq:two-bridge}
\khrn (K)(q,t) = (-t)^{-\sigma(K)/2} \pred(K)(q^Nt, iqt^{1/2}) \quad \mbox{for $N \geq 3$.}
\end{equation}
Most of
Conjecture~\ref{Conj:two-bridge} has been proved in
\cite{JakeUnwritten}, where it is shown that (\ref{eq:two-bridge})
holds for all $N \geq 5$.  The proof uses only elementary properties of
Khovanov and Rozansky's original definition, in particular the skein
exact sequence.  The approach has difficulties for $N = 3$ or $4$, and
this portion of Conjecture~\ref{Conj:two-bridge} remains open. All
knots with fewer than \(8\) crossings are two-bridge. Their
superpolynomials (assuming the conjecture) are shown in
Table~\ref{Tab:SPol}.

It is well known \cite{OSAlt}, \cite{ESL1} that alternating knots are
thin with respect to both \(\HFK\) and \(\hkr_2\). However, the
analogous statement for \(\hk\)--thinness cannot be true. To see why,
we introduce the notion of a knot having an {\it alternating} HOMFLY
polynomial. We say that \(P(K)\) is alternating if the sign of the
coefficient of \(a^{2i}q^{2j}\) is \(\pm (-1)^{j}\), where the factor of
\(\pm\) is the same for all coefficients. It is not difficult to see
that if \(K\) is \(\hk\)--thin, then \(P(K)\) is alternating.  On the
other hand, there are examples of alternating knots whose HOMFLY
polynomials are not alternating, the smallest being \(11^a_{263}\)
(numbering from {\it Knotscape} \cite{Knotscape}).

\begin{figure}
\centering\includegraphics{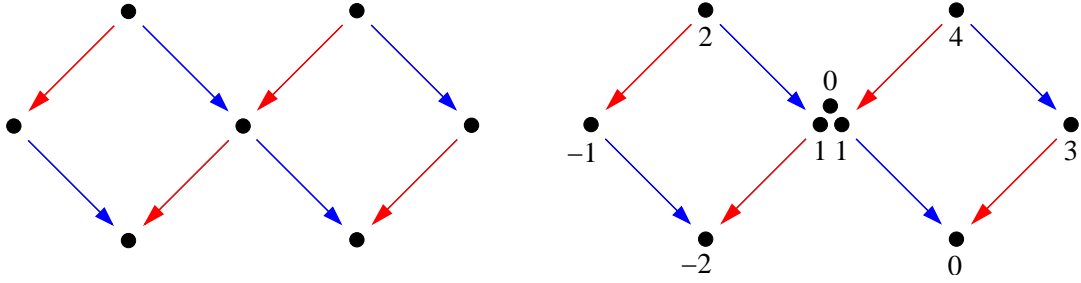}
\caption{\label{Fig:942} Two possible dot diagrams for the knot
\(9_{42}\). The left-hand diagram assumes that \(\hk(9_{42})\) is thin and
arrives at a contradiction: \(d_1^2 \neq 0\). The right-hand diagram
corrects this problem by introducing a pair of additional generators.}
\end{figure}

Conversely, knots with alternating HOMFLY polynomials need not be
\(\hk\)--thin. The knot \(9_{42}\) (numbering from Rolfsen \cite{Rolfsen}) is
a good example of this phenomenon. It has HOMFLY polynomial
\begin{equation}
P(9_{42}) = a^{-2}q^{-2} + a^{-2}q^{2} - q^{-4} - 1 - q^4 + a^2q^{-2} +
a^2q^2
\end{equation}
which is certainly alternating. If we assume \(\hk(9_{42})\) is thin
and try to endow it with differentials satisfying
Conjecture~\ref{Conj:Diffs}, however, we arrive at a contradiction.
The requirement that \(d_1\) and \(d_{-1}\) have one-dimensional
homology and anticommute with each other quickly leads to the dot
diagram shown on the left hand side of Figure~\ref{Fig:942}.  However,
in that diagram both \(d_1\) and \(d_{-1}\) do not square to zero. The
problem is resolved by postulating the presence of an additional two
generators in at the center of the diagram, as shown on the right-hand
side of Figure~\ref{Fig:942}. The resulting diagram correctly predicts
\(\HFK (9_{42})\) and \(\hkr_2 (9_{42})\).

\subsection{Thick knots}

Some knots are easily identified as being \(\hk\)--thick. In
particular, if a knot is thick with respect to either \(\HFK\) or
\(\hkr_2\), it is necessarily \(\hk\)--thick as well. The knots with
fewer than \(11\) crossings which fit this criterion are 

\begin{equation}
8_{19},
9_{42},
10_{124},10_{128},10_{132},10_{136},10_{139},
10_{145},10_{152},10_{153},10_{154},10_{161}.
\end{equation} 
We have already described the first two of these in
Figures~\ref{Fig:T34} and \ref{Fig:942}. In
Section~\ref{sec-dot-diagrams}, we give dot diagrams illustrating what
we believe are the superpolynomials of the \(10\)-crossing knots in
the list above.  For most of these knots, our reasons for asserting
that this is the superpolynomial are purely internal: it seems
difficult to produce another diagram satisfying all the hypotheses of
Conjecture~\ref{Conj:Diffs}. In addition, there are skein theoretic
arguments which support our calculations for \(8_{19}\) and
\(10_{128}\), although these currently fall short of a complete proof.
In both of these cases, the skein theoretic calculation gave the
answer we had previously guessed based on our conjecture, and we view
this as at least some evidence that our calculations are on the right
track.
 
%
%
%
%
%

The interesting examples provided by these thick knots allow us
to probe the rich structure of the triply graded theory.
Even  the simple thick knots we considered
exhibit some very different types of behavior. Some thick knots, like
\(9_{42}, 10_{132},10_{136}\) and \(10_{145}\) have ``invisible''
generators which cannot be seen from the HOMFLY polynomial. Others, like
\(8_{19}\) and \(10_{124}\) have no invisible generators, but have
nontrivial \(d_{-2},d_{0}\), and \(d_2\). Many exhibit both
features. There are cases, like \(10_{145}\), where the 
gradings in the superpolynomial
are such that \(d_2\) might conceivably be nontrivial, but the
requirement that the differentials anticommute prohibits it.

Although the sample of knots we consider here is admittedly small, a
number of interesting patterns may be observed from it. The rest of
this section is devoted to describing a few of these.

\subsection{Dimension of $\HFK$ and $\hkr_2$}\label{Subsec:ranks}

\begin{table}\begin{center}
\begin{tabular}{|c|l|}
\hline
\rule{0pt}{5mm}
Knot & $\Pred$ \\[3pt]
\hline
\hline
\rule{0pt}{5mm}
$3_1$ & $a^2 q^{-2} + a^2 q^2 t^2 + a^4 t^3$ \\[3pt]
\hline
\rule{0pt}{5mm}
$4_1$ & $a^{-2} t^{-2} + q^{-2} t^{-1} + 1 + q^2 t + a^2 t^2$ \\[3pt]
\hline
\rule{0pt}{5mm}
$5_1$ & $a^4 q^{-4} + a^4 t^2 + a^6 q^{-2} t^3 + a^4 q^4 t^4
+ a^6 q^2 t^5$ \\[3pt]
\hline
\rule{0pt}{5mm}
$5_2$ & $a^2 q^{-2} + a^2 t + a^2 q^2 t^2 + a^4 q^{-2} t^2
+ a^4 t^3 + a^4 q^2 t^4 + a^6 t^5$ \\[3pt]
\hline
\rule{0pt}{5mm}
$6_1$ & $a^{-2} t^{-2} + q^{-2} t^{-1} + 2 + q^2 t + a^2 q^{-2} t
+ a^2 t^2 + a^2 q^2 t^3 + a^4 t^4$ \\[3pt]
\hline
\rule{0pt}{5mm}
$6_2$ & $q^{-2} t^{-2} + a^2 q^{-4} t^{-1} + a^2 q^{-2} + q^2
+ 2 t a^2 + a^4 q^{-2} t^2 + a^2 q^2 t^2$ \\[3pt]
\rule{0pt}{5mm}
& $+ a^2 q^4 t^3 + a^4 t^3 + a^4 q^2 t^4$ \\[3pt]
\hline
\rule{0pt}{5mm}
$6_3$ & $a^{-2} q^{-2} t^{-3} + a^{-2} t^{-2} + q^{-4} t^{-2}
+ q^{-2} t^{-1} + a^{-2} q^2 t^{-1} + 3$ \\[3pt]
\rule{0pt}{5mm}
& $+ a^2 q^{-2} t + q^2 t + a^2 t^2 + q^4 t^2 + a^2 q^2 t^3$ \\[3pt]
\hline
\rule{0pt}{5mm}
$7_1$ & $a^6 q^{-6} + a^6 q^{-2} t^2 + a^8 q^{-4} t^3
+ a^6 q^2 t^4 + a^8 t^5 + a^6 q^6 t^6 + a^8 q^4 t^7$ \\[3pt]
\hline
\rule{0pt}{5mm}
$7_2$ & $a^2 q^{-2} + a^2 t + a^4 q^{-2} t^2 + a^2 q^2 t^2
+ 2 a^4 t^3 + a^6 q^{-2} t^4 + a^4 q^2 t^4$ \\[3pt]
\rule{0pt}{5mm}
& $+ a^6 t^5 + a^6 q^2 t^6 + a^8 t^7$ \\[3pt]
\hline
\rule{0pt}{5mm}
$7_3$ & $a^{-4} q^4 + a^{-8} q^{-2} t^{-7} + a^{-6} q^{-4} t^{-6}
+ a^{-6} q^{-2} t^{-5} + a^{-8} q^2 t^{-5} + 2 a^{-6} t^{-4}$ \\[3pt]
\rule{0pt}{5mm}
& $+ a^{-4} q^{-4} t^{-4} + a^{-6} q^{-2} t^{-3} + a^{-6} q^2 t^{-3}
+ a^{-4} t^{-2} + a^{-6} q^4 t^{-2} + a^{-4} q^2 t^{-1}$ \\[3pt]
\hline
\rule{0pt}{5mm}
$7_4$ & $a^{-2} q^2 + a^{-8} t^{-7} + a^{-6} q^{-2} t^{-6}
+ 2 a^{-6} t^{-5} + 2 a^{-4} q^{-2} t^{-4} + a^{-6} q^2 t^{-4}$ \\[3pt]
\rule{0pt}{5mm}
& $+ 2 a^{-4} t^{-3} + a^{-2} q^{-2} t^{-2} + 2 a^{-4} q^2 t^{-2}
+ 2 a^{-2} t^{-1}$ \\[3pt]
\hline
\rule{0pt}{5mm}
$7_5$ & $a^4 q^{-4} + a^4 q^{-2} t + 2 a^4 t^2 + a^6 q^{-4} t^2
+ 2 a^6 q^{-2} t^3 + a^4 q^2 t^3 + 2 a^6 t^4$ \\[3pt]
\rule{0pt}{5mm}
& $+ a^4 q^4 t^4
+ a^8 q^{-2} t^5 + 2 a^6 q^2 t^5 + a^8 t^6 + a^6 q^4 t^6
+ a^8 q^2 t^7$ \\[3pt]
\hline
\rule{0pt}{5mm}
$7_6$ & $2 a^2 q^{-2} + q^2 + q^{-2} t^{-2} + t^{-1}
+ a^2 q^{-4} t^{-1} + 3 a^2 t + 2 a^4 q^{-2} t^2$ \\[3pt]
\rule{0pt}{5mm}
& $+ 2 a^2 q^2 t^2 + 2 a^4 t^3 + a^2 q^4 t^3 + 2 a^4 q^2 t^4 + a^6 t^5$ \\[3pt]
\hline
\rule{0pt}{5mm}
$7_7$ & $a^{-4} t^{-4} + 2 a^{-2} q^{-2} t^{-3}
+ 2 a^{-2} t^{-2} + q^{-4} t^{-2} + 2 q^{-2} t^{-1}
+ 2 a^{-2} q^2 t^{-1}$ \\[3pt]
\rule{0pt}{5mm}
& $ + 4 + a^2 q^{-2} t + 2 q^2 t
+ 2 a^2 t^2 + q^4 t^2 + a^2 q^2 t^3$ \\[3pt]
\hline
\end{tabular}\end{center}
\caption{\label{Tab:SPol} Reduced superpolynomial for prime knots
with up to 8 crossings.}
\end{table}

It is an interesting and rather puzzling fact that the knot Floer
homology and $sl(2)$ Khovanov homology of a given knot often have the
same dimension \cite{KFC}. Indeed, explaining this was one of our
motivations for considering a triply graded theory. At first glance,
however, the triply graded theory we have described does not seem to
help all that much.  One case where it does provide insight is for
those knots where \(d_2\) and \(d_0\) both vanish (thin knots, but
also some thick examples such as \(9_{42}\)).  In this case, the
correspondence is obvious: the dimensions of $\HFK$ and $\hkr_2$ are
both equal to that of \(\hk\).  However, there are many knots where
\(d_2\) and \(d_0\) are nontrivial but the two dimensions still agree.
To consider an extreme example, our proposal for \(\hk(10_{128})\) has
dimension \(27\), while the dimension of \(\HFK\) and \(\hkr_2\) are
both 13.

The fact that the correspondence still holds in such cases suggests
that we should look for an explanation of why the part of \(\hk\)
killed by \(d_2\) should have the same dimension as the part killed by
\(d_0\). Examining the diagrams in Section~\ref{sec-dot-diagrams}, a
rather striking pattern comes to light: for knots with \(S \geq 0 \),
any dot that has a nonzero image under one of \(d_2,d_0\), and
\(d_{-2}\) must have a nonzero image under the other two as well! (For
\(S<0\), the requirement is reversed: any generator that is in the
image of one differential is in the image of the other two as well.)
Although we don't have any explanation for this phenomenon, it seems
clear that if we understood it, we would be well on the way to
understanding why \(\HFK\) and \(\hkr_2\) have the same dimension for
so many knots.
 
\subsection{Braid index and estimates on $S$}\label{Subsec:four-ball}

It is well known that the minimum braid index of a knot is bounded by
the difference between the maximum and minimum exponents of \(a\) in
its HOMFLY polynomial. The same principle applies to the
superpolynomial. More generally, we have 
\begin{prop}
\label{Prop:Morton}
Let \(a_{\max}(\Pred(K))\) and \(a_{\min}(\Pred(K))\) be the maximum and minimum powers of
\(a\) appearing in \(\Pred (K)\). Then for any planar diagram \(D\) of
\(K\), 
\begin{equation}
w(D)-c(D) + 1 \leq a_{\min}(\Pred(K)) \leq a_{\max}(\Pred(K)) 
\leq w(D) + c(D) -1
\end{equation}
where \(w(D)\) is the writhe of \(D\) and \(C(D)\) the number
of components in its oriented resolution. 
\end{prop}

The analog of this theorem for the HOMFLY polynomial was proved by
Morton in \cite{Morton}.  As we now describe, Morton's argument
carries through to the setting of superpolynomials.  Since we don't
have a definition of $\Pred(K)$, this statement can be taken in two
ways.  The first is that, like the $sl(N)$ homology, the triply graded
theory should satisfy a skein exact triangle.  Morton's proof is
purely skein-theoretic, and it is not hard to see it carries over to
any theory that has a skein exact triangle.  The other point of view
is that this is a limiting statement about the $sl(N)$ homology as $N
\to \infty$.   In particular, using the skein exact triangle one can show 
\begin{equation}\label{eq:approx-morton}
  N \left( w(D) - c(D) + 1 \right) - E \leq q_{\min}(\khr_N(K)) \leq q_{\max}(\khr_N(K)) \leq N \left( w(D) + c(D)  -1 \right) + E.
\end{equation}
where $\abs{E}$ is uniformly bounded independent of $N$.  Provided
Conjecture~\ref{khredconj} holds, we have $\lim_{N \to \infty} (1/N)
q_{\min}(\khr_N(K)) = a_{\min}(\Pred(K))$ and similarly for
$a_{\max}(\Pred(K))$.  The proposition then follows by taking the limit
of \eqref{eq:approx-morton} as $N \to \infty$.

In the same paper, Morton asked whether there might be a connection
between \(a_{\min}(P(K))\) and the bound on the genus of a knot
provided by Bennequin's inequality. Since Bennequin's inequality
actually provides a lower bound for the four-ball genus $g_*$ of \(K\)
\cite{Rudolph}, one might ask whether the same is true for
\(a_{\min}(P(K))\):
\begin{equation}\label{eq-4-genus-bound}
2g_*(K) \stackrel{??}{\geq} a_{\min}(P(K)).
\end{equation} 
Although it is true in many examples, this inequality is false in
general.  For knots with fewer than \(11\) crossings, the knot
\(K=10_{132}\) is the only counterexample; there \(g_*(K)=1\), but
\(a_{\min}(P(K))=4\). A brief inspection of the proposed dot diagram
for \(10_{132}\) in Section~\ref{sec-dot-diagrams} suggests an
explanation for what has gone wrong: \( a_{\min}(\Pred(K))=2\), but
the terms with lowest degree in \(a\) are not visible in the HOMFLY
polynomial.

If we replace \(a_{\min}(P(K))\) by \(a_{\min}(\Pred(K))\) in
\eqref{eq-4-genus-bound}, we expect that the resulting inequality will
be true. Indeed, it is clear from the definition that
\(a_{\min}(\Pred(K)) \leq S(K) \leq a_{\max}(\Pred(K))\). If \(S(K)\)
provides a lower bound for the four-ball genus of \(K\) (which seems
quite likely), \( a_{\min}(\Pred(K))\) will do so as well.  Continuing
in this vein, we can combine Proposition~\ref{Prop:Morton} with the
previous inequality to obtain the following estimate for \(S\):
\begin{equation}
w(D)-c(D) + 1 \leq S(K) \leq w(D) +c(D) - 1
\end{equation}
where \(D\) is any planar diagram of \(K\).  Zolt{\'a}n Szab{\'o} pointed
out to us that using the work of Livingston \cite{Livingston}, it is not
difficult to see that \(s\) and \(\tau\) satisfy similar estimates. We
sketch the proof of this fact for \(\tau\); the argument for \(s\) is
the same.

Suppose \(K\) has a planar diagram \(D\), and let
\(n_{\pm}(K)\) denote the number of positive and negative
crossings. If we change all the negative crossings to positive, we
obtain a new knot \(K^+\), and \cite{Livingston} and \cite{Rudolphii}
tell us that
\begin{equation}
2 \tau(K^+) = n_{+}(D) + n_{-}(D) - c(D)+1.
\end{equation}
To get back to \(K\), we must change \(n_-(D)\) crossings from
positive to negative, which can lower \(\tau\) by at most
\(n_{-}(D)\). Thus
\begin{equation}
2 \tau (K) \geq  {n_{+}(D) - n_{-}(D) - c(D)+1} = {w(D) - c(D)+1}.
\end{equation}
Similarly, changing all of \(D\)'s
positive crossings to negative, we see that 
\begin{equation}
2 \tau (K) \leq {w(D) + c(D)+1}. 
\end{equation}


\subsection{ $d_1$ and the unreduced homology}

Although we have focused on reduced homology, we expect that our work
also has relations with the unreduced theory. In general, the
unreduced homology \(\hkun_N(K)\) is related to \(\hkr_N (K)\) by a
spectral sequence which has \(E_1\) term equal to \(\hkr_N(K) \otimes
\Q[X]/(X^N)\). When \(N=2\), the differential in this spectral
sequence seems to be related to the Lee/Turner differential on
\(\hkr_2\). For example, if \(K\) is thin, the presence of the
Lee/Turner differential implies that
\begin{equation}\label{one-step-pairing}
\khr_2(K) = q^{s(K)} + (1+q^2t) \khr_2'(K),
\end{equation}
where \(\khr_2'(K) \) is a polynomial with positive coefficients. The unreduced
homology can also be expressed in terms of \(\khr_2'\): 
\begin{equation}
\label{Eq:Kh2un}
\kh_2 (K) = (q+q^{-1}) q^{s(K)} + (q^{-1}+q^3t)\khr_2'(K).
\end{equation}
This suggests that the differential on the \(E_1\) term of 
 the spectral sequence is determined by the relation
\(d_{E_1}(a) = X d_{1*}(a)\), where \(d_{1*}\) denotes the
Lee/Turner differential.

The analog for the superpolynomial is that for \emph{any} knot $K$ we
have:
\begin{equation}
\label{patternb}
\Pred(K)  = \Big( \frac{a}{q} \Big)^{S(K)} + 
(1+ t a^2 q^{-2} ) Q_+(a,q,t)
\end{equation}
where $Q_{+}(a,q,t)$ is a polynomial with positive coefficients.  This
follows immediately from the existence of the canceling differential
$d_1$ given by Conjecture~\ref{Conj:Diffs}.  (The reason that the
standard canceling differential on $\hkr_2$ does not always force
(\ref{one-step-pairing}) is that, unlike $d_1$ on $\hk_*$, it is not
necessarily homogeneous in its behavior with respect to the grading.)
When \(K\) is thin, we expect that the differential in the spectral
sequence will again be determined by \(d_{1*}\): \(d_{E_{N-1}} (a) =
X^{N-1} d_{1*}(a)\).  This suggests the following analog of
\eqref{Eq:Kh2un}:
\begin{equation}
\kh_N(K) = q^{(N-1)S(K)}\left( \frac{q^N-q^{-N}}{q-q^{-1}}\right)
 +(q^{-1} + q^{2N-1}t) \left( \frac{q^{N-1}-q^{-N+1}}{q-q^{-1}} \right)
Q_+(a=q^{N},q,t).
\end{equation}
Expressing this in terms of the unreduced superpolynomial, we get
\begin{equation}
\label{patterna}
\Pun(K) = (a - a^{-1}) \Big( \frac{a}{q} \Big)^{S(K)}  + 
(q^{-1}+a^2q^{-1}t)(aq^{-1}-a^{-1}q)Q_{+}(a,q,t). 
\end{equation}
Let us illustrate the structure of the unreduced superpolynomial
with the following example.

\begin{example:figure8}
Since the figure-eight knot $4_1$ is $\hk$-thin, its reduced superpolynomial
is easy to determine. The result is presented in Table \ref{Tab:SPol}.
It has the expected structure \eqref{patternb} with $S(4_1)=0$ and
\begin{equation}\label{figeightpredprime}
\Pred'(4_1) = {1 \over a^2 t^2} + q^2 t
\end{equation}
Substituting this into \eqref{patterna}, we find the unreduced
superpolynomial for the figure-eight knot:
\begin{equation}\label{punfigeight}
\Pun (4_1) = a - a^{-1}
+ (q^{-1}+a^2q^{-1}t)(aq^{-1}-a^{-1}q)(a^{-2} t^{-2} + q^2 t)
\end{equation}
It is easy to check that specializing to $t=-1$ and $a=q^2$
we reproduce , respectively, the correct expressions for
the unnormalized HOMFLY polynomial and the $sl(2)$ Khovanov homology.
Moreover, substituting \eqref{punfigeight} into \eqref{khnviap},
we obtain the following prediction for the unreduced $sl(N)$ homology:
\begin{equation}
\label{slnfigeight}
\kh_N(4_1) = \sum_{i=0}^{N-1} q^{2i-N+1} + 
(1+q^{2N}t)(q^{-2N} t^{-2} + q^2 t) \sum_{i=0}^{N-2} q^{2i-N+1}
\end{equation}
\end{example:figure8}


\subsection{$d_{-1}$ and three-step pairings}
\label{Subsec:d-1}

As discussed in Section~\ref{Subsec:CDiffs},
Conjecture~\ref{Conj:Diffs} requires that \(\hk\) admit two distinct
canceling differentials: \(d_1\) and \(d_{-1}\). This implies that
\(\hkr_N\) should admit a second canceling differential as well.  We
end this section by describing some empirical evidence which supports
the idea that \(\hkr_2\) admits an additional canceling differential.
 
To begin with, we show that the unique term that is not canceled by
$d_{-1}$ has grading \((aqt)^{S(K)}\).  This is because
$d_{-1}$ is interchanged with $d_1$ by the symmetry $\phi$ --- the
uncanceled term for $d_1$ is $a^{S(K)}q^{-S(K)}t^0$ which is taken to
\(a^{S(K)}q^{S(K)}t^n\) by $\phi$, and $n$ can then be computed by
using that $\phi$ preserves the $\delta$-grading. We thus have the
following analog of \eqref{patternb}:
\begin{equation}
\label{patternbb}
\Pred(K) = (aqt)^{S(K)} + (1+a^2q^2t^3) Q_-(a,q,t)
\end{equation}
where \(Q_-(a,q,t)\) is a polynomial with positive coefficients. 

If \(K\) is $\hk$-thin, we can substitute \(a=q^2\) to obtain the following
 prediction
 for the \(sl(2)\) Khovanov homology of \(K\):
\begin{equation}\label{patternbj}
\khr_2(K) = (q^3t)^{S(K)} + (1+q^6t^3) Q_-(a=q^2,q,t).
\end{equation}
Independent of this, given a $\hkr_2$-thin knot \(K\), we have
\begin{equation}
\khr_2(K) = (-t)^{-S(K)/2} J(K)(q^2 = -q^2t)
\end{equation}
where $J(K)$ is the Jones polynomial $P(K)(a=q^2, q)$.  Combining this
with the fact that \(J(K)(q^2)-1\) is divisible by \(1-q^6\) (see {\it
  e.g.} Proposition 12.5 of \cite{Jonestorus}), it is not difficult to
see that \eqref{patternbj} holds for some polynomial \(Q_-(q,t)\).
It's not clear that this polynomial should have positive
coefficients, as predicted by \eqref{patternbb}, but for thin knots with
fewer than \(12\) crossings, we have checked that this is the case.
More generally, we make the following

\begin{definition}
We say a knot \(K\) has a \emph{three-step pairing} on \(\khr_2\)
 if for some \(m,n \in \Z\), we have 
\begin{equation}
\khr_2(K) = q^mt^n + (1+q^6t^3)Q_-(q,t)
\end{equation} 
 where \(Q_-\) is a polynomial with positive coefficients.
\end{definition}

A knot which admits a three-step pairing has an obvious candidate for
the canceling differential induced by \(d_{-1}\), though of course a
canceling differential need not force a three-step pairing.  Such
knots are surprisingly common. In addition to the thin knots mentioned
above, we checked some $5{,}000$ knots with fewer than
than 16 crossings which happen to be $(1,1)$ knots and found that all
of them had three-step pairings.  A number of these knots are complicated
enough that they do not satisfy \eqref{one-step-pairing}, which makes
this all the more remarkable.


\section{Torus knots}\label{sec-torus-1}

Let $T_{n,m}$ be a torus knot of type $(n,m)$,
where $n$ and $m$ are relatively prime integers, $n<m$.
In this section, we propose an explicit expression for
the superpolynomial for all torus knots of type $(2,m)$ and $(3,m)$,
and discuss its structure for general torus knots $T_{n,m}$.
We consider reduction to the $sl(N)$ knot homology
and to the knot Floer homology, and show that our
predictions are consistent with the known results.
The differentials $d_N$ play an important role in this discussion.

Let us begin by recalling the expression
for the HOMFLY polynomial of a torus knot $T_{n,m}$.

\subsection{HOMFLY polynomial}

The explicit expression for $\pred (T_{n,m})$
was found by Jones \cite{Jonestorus}:
\begin{equation}
\label{homflyfortnm}
\pred (T_{n,m}) = {a^{m(n-1)} [1]_q \over [n]_q}
\sum_{\b=0}^{n-1} (-1)^{n-1-\b}
{q^{-m(2\b-n+1)} \over [\b]_q ! [n-1-\b]_q !}
\prod_{{j=\b-n+1 \atop j \neq 0}}^{\b} \big(q^{j} a - q^{-j} a^{-1} \big)
\end{equation}
where $[n]_q = q^n - q^{-n}$ is the ``quantum dimension'' of $n$
written in a slightly unconventional normalization, and
\begin{equation}
\label{qfactorial}
[n]_q ! = [n]_q [n-1]_q \ldots [1]_q
\quad \mbox{with} \quad
[0]_q ! = 1.
\end{equation}
One can manipulate the expression \eqref{homflyfortnm} into
the following form, which will be useful to us below,
\begin{equation}
\label{homflyfortnmnew}
\pred (T_{n,m})
= (aq)^{(n-1)(m-1)} {1 - q^{-2} \over 1 - q^{-2n}}
\sum_{\b=0}^{n-1} q^{-2m\b}
\left( \prod_{i=1}^{\b} {a^2 q^{2i}-1 \over q^{2i}-1} \right)
\left( \prod_{j=1}^{n-1-\b} {a^2 - q^{2j} \over 1-q^{2j}} \right).
\end{equation}

Assuming that all the terms in the superpolynomial $\Pred (T_{n,m})$
are ``visible'' in the HOMFLY polynomial, one might hope to obtain
$\Pred (T_{n,m})$ by inserting powers of $(-t)$ in the expression
for $\pred (T_{n,m})$. In order to do this, it is convenient
to simplify \eqref{homflyfortnm} further and write it as
a sum of terms without denominators.
For example, for $n=2$ and $m=2k+1$, we find
\begin{eqnarray}
\label{homflyforntwo}
\pred (T_{2,2k+1}) & = & {a^{2k+1} \over (q^2 - q^{-2})}
\Big[ - a (q^{2k} - q^{-2k})
+ a^{-1} \big( q^{2k+2} - q^{-2k-2} \big) \Big] = \cr
& = & - a^{2k+2} \sum_{i=1}^k q^{4i-2k-2}
+ a^{2k} \sum_{i=0}^k q^{4i-2k}
\end{eqnarray}
where in the first line we combined the terms with the same
power of $a$. Similarly, for $(3,m)$ torus knots, we find
\begin{equation}
\label{homflyfornthreea}
\pred (T_{3,3k+1}) = 
a^{6k} \sum_{j=0}^{k} \sum_{i=0}^{3j} q^{6j-4i}
- a^{6k+2}\sum_{j=1}^{k} \sum_{i=0}^{6j-2}  q^{6j-2i-2}
+ a^{6k+4} \sum_{j=0}^{k-1} \sum_{i=0}^{3j}  q^{6j-4i}
\end{equation}
and
\begin{equation}
\label{homflyfornthreeb}
\pred (T_{3,3k+2}) = 
a^{6k+2} \sum_{j=0}^{k} \sum_{i=0}^{3j+1}  q^{6j-4i+2}
- a^{6k+4} \sum_{j=0}^{k} \sum_{i=0}^{6j}  q^{6j-2i}
+ a^{6k+6} \sum_{j=0}^{k-1} \sum_{i=0}^{3j+1}  q^{6j-4i+2}
\end{equation}
In general, $\pred (T_{n,m})$ has the following structure,
which follows directly from \eqref{homflyfortnmnew},
\begin{equation}
\label{ptnmstructure}
\pred (T_{n,m}) = \sum_{J=0}^{n-1} a^{(m-1)(n-1)+2J} \pred^{(J)} (q)
\end{equation}
where each $\pred^{(J)} \in \Z[q,q^{-1}]$
can be written in terms of $n-1$ repeated sums,
{\it cf.} \eqref{homflyforntwo} -- \eqref{homflyfornthreeb}.


\subsection{The structure of the superpolynomial}

We wish to find an explicit form of the superpolynomial
for torus knots $T_{n,m}$, which has all the right
properties to be the Poincar\'e polynomial of 
the triply graded homology theory $\hk$.
Before we proceed to a more detailed analysis,
let us make a few general remarks about the expected
structure of the superpolynomial for torus knots $T_{n,m}$.
Simple examples of torus knots of type $(2,m)$
and $(3,m)$ already appeared in Sections \ref{sec-differentials}
and \ref{Sec:Examples}.
In these examples, all the terms in the reduced superpolynomial $\Pred (T_{n,m})$
are ``visible'' in the HOMFLY polynomial. We will assume
that this is also the case for more general torus knots.
In particular, this means that the superpolynomial $\Pred (T_{n,m})$
has the structure similar to \eqref{ptnmstructure},
\begin{equation}
\label{superptnm}
\Pred (T_{n,m}) = \sum_{J=0}^{n-1} a^{(m-1)(n-1)+2J} \Pred^{(J)} (q,t)
\end{equation}
where
$$
\Pred^{(J)} \in \Z_{\geq 0} [q,q^{-1},t]
$$
Notice that only non-negative powers of $t$ appear in $\Pred^{(J)}(q,t)$.
Moreover, the examples of $T_{2,m}$ and $T_{3,m}$ torus knots
studied below suggest that only even (resp.~odd) powers of $t$
appear in $\Pred^{(J)}(q,t)$ for even (resp.~odd) values of $J$,
and the maximal degree of $t$ does not exceed $(m-1)(n-1)+J$.

The structure of the superpolynomial $\Pred (T_{n,m})$ should be
also consistent with the action of the differentials $d_1$ and $d_{-1}$.
In particular, it should be consistent
with \eqref{patternb} and \eqref{patternbb}:
\begin{equation}
\label{torussfirst}
\Pred (T_{n,m}) = a^S q^{-S} + (a^2 q^{-2} t + 1) Q_{+} (a,q,t)
\end{equation}
\begin{equation}
\label{torusslast}
\Pred (T_{n,m}) = (aqt)^S + (a^2 q^2 t^3 + 1) Q_{-} (a,q,t)
\end{equation}
where, for a torus knot $T_{n,m}$,
\begin{equation}
\label{sfortnm}
S(T_{n,m}) = (n-1)(m-1)
\end{equation}
and $Q_{\pm} \in \Z_{\geq 0} [a,q,t]$.
Similarly, the unreduced superpolynomial should have
the structure, {\it cf.} \eqref{patterna}:
\begin{equation}
\label{torusfirstun}
\Pun (T_{n,m}) = \left( {a \over q} \right)^S (a - a^{-1})
+ (a^{-1} + at) \Pun' (a, q, t)
\end{equation}
where $\Pun' \in \Z [a^{\pm 1},q^{\pm 1},t^{\pm 1}]$.

We believe that, for any torus knot $T_{n,m}$,
there exists an explicit expression for
the superpolynomial with all the required properties.
We were able to find such an expression for all torus
knots of type $(2,m)$ and $(3,m)$, and to obtain some
partial results for arbitrary torus knots $T_{n,m}$.


\subsection{Torus knots $T_{2,2k+1}$}
\label{Subsec:T(2,m)}

The \((2,2k+1)\) torus knots are in many respects the simplest of all
knots.  There are several different ways to determine their
superpolynomials (reduced and unreduced), all of which lead to the
same result.  One reason for this --- which was already used for
simple examples of $(2,2k+1)$ torus knots in \cite{GSV} and in
Sections \ref{sec-differentials} and \ref{Sec:Examples} here --- is
that all the terms in the $sl(2)$ homology of $T_{2,2k+1}$ are
``visible'' in the HOMFLY polynomial.  In particular, for torus knots
of type $(2,2k+1)$, the Conjectures \ref{khnconj} and \ref{khredconj}
hold for all values of $N \geq 2$.  This nice property can be used to
determine the superpolynomial of $T_{2,2k+1}$ either by combining the
information about the HOMFLY polynomial and the $sl(2)$ homology, or
by comparing the $sl(2)$ and $sl(3)$ knot homologies, or in some other
way.

For example, the HOMFLY polynomial of \(T_{2,2k+1}\) is given by 
\eqref{homflyforntwo}:
\begin{equation}
\pred (T_{2,2k+1})  =  - a^{2k+2} \sum_{i=1}^k q^{4i-2k-2}
+ a^{2k} \sum_{i=0}^k q^{4i-2k}
\end{equation}
while the \(sl(2)\) Khovanov homology is
\begin{equation}
\khr_2(T_{2,2k+1}) = q^{2k} t^0 + q^{2k+4}t^2 + q^{2k+6}t^3+\ldots q^{6k+2}t^{2k+1}.
\end{equation}
If we substitute \(a=q^2\) and compare terms, it is easy to guess the 
following formula:

\begin{prop}\label{Prop:ttwomreed}
The reduced superpolynomial $\Pred (T_{2,2k+1})$
has the form \eqref{superptnm}:
\begin{equation}
\label{predformtwo}
\Pred (T_{2,2k+1}) =
a^{2k} \Pred^{(0)} + a^{2k+2} \Pred^{(1)}
\end{equation}
where 
\begin{equation}
\label{ppformtwo}
\Pred^{(0)} = \sum_{i=0}^k q^{4i-2k} t^{2i} \quad \mbox{and} \quad
\Pred^{(1)} = \sum_{i=1}^k q^{4i-2k-2} t^{2i+1}.
\end{equation}
\end{prop}
Of course, \(T_{2,2k+1}\) is a two bridge knot, so  
 a particular case  of Conjecture~\ref{Conj:two-bridge}. This is a very
 useful family of examples to have in mind, however, so it is worth
 considering them in greater detail. 
Note that we have stated the formula above as a proposition. As
usual, this is to be interpreted as a statement about \(\khrn\) for
\(N \gg 0\). Its proof follows immediately from the proof of
\ref{Conj:two-bridge} given in \cite{JakeUnwritten}.

\begin{figure}
\centering\includegraphics[scale=1.2]{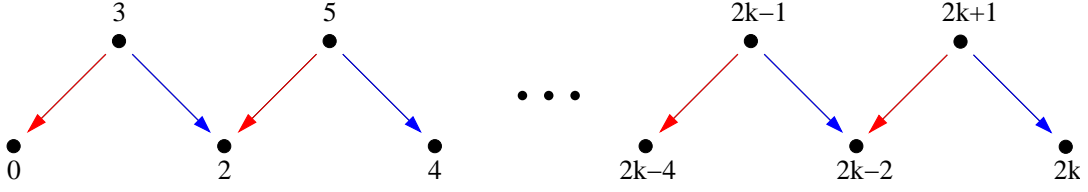}
\caption{\label{Fig:T22k1}
Dot diagram for the superpolynomial of \(T_{2,2k+1}\).}
\end{figure}

Let us check that \(\hk(T_{2,2k+1})\) satisfies the conditions of
Conjecture~\ref{Conj:Diffs}. First, observe that \(\hk(T_{2,2k+1})\)
is thin --- all generators have \(\delta\)-grading \(-k\). For \(i\neq
0\),  \(d_i\) lowers the \(\delta\)-grading by \(|i|\), while \(d_0\)
lowers the \(\delta\)-grading by 1. 
Thus \(d_1\) and \(d_{-1}\) must be the only nontrivial differentials.
Their action is illustrated in Figure~\ref{Fig:T22k1}.
 From the figure, it is obvious that the symmetry property holds. 
Finally, if we substitute  $a=1/t$, the reduced superpolynomial 
specializes to $\hfk (T_{2,2k+1})$:
\begin{equation}
\label{hfkfortntwo}
\hfk (T_{2,2k+1}) = q^{-2k} t^{-2k}
+ q^{-2k} t^{-2k} (1 + q^{-2} t^{-1}) \sum_{i=1}^k q^{4i} t^{2i}
\end{equation}

We remark that the vanishing of \(d_N\) for \(N\neq 1, -1\) is really quite special. 
As we shall see in the next section, the situation is qualitatively
different for torus knots $T_{n,m}$ with $n>2$,
where any differential $d_N$ can potentially be
non-trivial for fixed value of $n$ and sufficiently large $m$.


Now, let us turn to the unreduced superpolynomial of $T_{2,2k+1}$.
The unnormalized HOMFLY polynomial of $T_{2,2k+1}$ can be
easily obtained from \eqref{homflyforntwo} by multiplying
it with $\pun(\mathrm{unknot}) = (a-a^{-1})/(q-q^{-1})$:
\begin{equation}
\label{unhomflyforntwo}
\pun (T_{2,2k+1}) = {1 \over q-q^{-1}} \left[
- a^{2k+3} \sum_{i=1}^k q^{4i-2k-2}
+ a^{2k+1} \sum_{i=0}^{2k} q^{2i-2k}
- a^{2k-1} \sum_{i=0}^k q^{4i-2k}
\right]
\end{equation}
On the other hand, the unreduced \(sl(2)\) homology of $T_{2,m}$
is known to be given by \cite{Khovanov} 
\begin{equation}
\label{khtwoforntwo}
\kh_2 (T_{2,2k+1}) = (q + q^{-1}) q^{2k}
+ \sum_{i=1}^k q^{4i+2k-1} t^{2i}
+ \sum_{i=1}^k q^{4i+2k+3} t^{2i+1}
\end{equation}
Now one can use the conjectured relation \eqref{khnviap}
to find the superpolynomial $\Pun (T_{2,2k+1})$.
Namely, multiplying both \eqref{unhomflyforntwo}
and \eqref{khtwoforntwo} by $(q-q^{-1})$ we obtain
two expressions, which are supposed to be specializations
of $\Pun (T_{2,2k+1})$ to $t=-1$ and $a=q^2$, respectively:
\begin{align}
& (q-q^{-1})\pun (T_{2,2k+1}) =
(a - a^{-1}) \left( {a \over q} \right)^{2k}
+ a^{2k} (aq^{-2} - a^3q^{-2}-a^{-1}+a) \sum_{i=1}^k q^{4i-2k} \\
& (q-q^{-1}) \kh_2 (T_{2,2k+1}) = (q^2 - q^{-2}) q^{2k}
+ (1 + q^4t - q^{-2} - q^2 t) \sum_{i=1}^k q^{4i+2k} t^{2i}
\end{align}
Matching the corresponding terms in
these two expressions, we arrive at the following formula, which is a special case of \eqref{patterna}:

\begin{prop}\label{Prop:ttwom}
For a torus knot $T_{2,2k+1}$,
the unreduced superpolynomial $\Pun (T_{2,2k+1})$ is given by
\begin{equation}
\label{punforttwom}
\Pun (T_{2,2k+1}) =
(a - a^{-1}) \left( {a \over q} \right)^{2k}
+ a^{2k}(a^2 q^{-2} - 1) (a^{-1} + at)
\sum_{i=1}^k q^{4i-2k} t^{2i}.
\end{equation}
\end{prop}

As a mathematical statement, this is to be interpreted in terms of Conjecture~\ref{khnconj}. In other words, it says that for \(N>1\), the 
 the $sl(N)$ knot homology of $T_{2,2k+1}$ is given by:
\begin{equation}
\label{khnfortntwo}
\kh_N (T_{2,2k+1}) = q^{(2k-1)(N-1)} \left[ \sum_{i=0}^{N-1} q^{2i}
+ (1 + q^{2N}t) \sum_{i=1}^k \sum_{j=0}^{N-2} q^{4i+2j} t^{2i} \right].
\end{equation}
Again, this formula can be confirmed by direct calculation. Perhaps the easiest approach is to start from Proposition~\ref{Prop:ttwomreed} and use the spectral sequence relating reduced and unreduced homology. All the differentials in this spectral sequence vanish for dimensional reasons except for \(d_{N-1}\), which is potentially nonzero on \(k\) different elements. To verify the nontriviality of \(d_{N-1}\), one can use Gornik's theorem \cite{Gornik} that there is a differential on \(\hkr_N(K)\) whose homology is supported in dimension zero. This cannot be the case unless all components of \(d_{N-1}\) which can be nonzero actually are nonzero.


\subsection{Torus knots $T_{3,m}$}

In this and the following section,
we consider torus knots of type $(3,m)$,
and we will mainly discuss the reduced theory.
We start by summarizing our prediction for
the superpolynomial of $T_{3,m}$:

\begin{conj}\label{Conj:tthreem}
For a torus knot $T_{3,m}$,
the reduced superpolynomial $\Pred (T_{3,m})$ has the form
\eqref{superptnm}:
\begin{equation}
\label{predforttrim}
\Pred (T_{3,m}) = a^{2m-2} \Pred^{(0)}
+ a^{2m} \Pred^{(1)} + a^{2m+2} \Pred^{(2)}
\end{equation}
where for $m=3k+1$
\begin{align}
\Pred^{(0)} &= \sum_{j=0}^{k} \sum_{i=0}^{3j}
q^{6j-4i} t^{4k+2j-2i}  \notag  \\
 \Pred^{(1)} &= \sum_{j=1}^{k} \sum_{i=0}^{6j-2}
q^{6j-2i-2} t^{4k+2j-2\lfloor i /2\rfloor+1} \quad \mbox{where $\lfloor x\rfloor$ denotes the integer part of $x$,} \label{pjforkone} \\ 
 \Pred^{(2)} &= \sum_{j=0}^{k-1} \sum_{i=0}^{3j}
q^{6j-4i} t^{4k+2j-2i+4} \notag \\
\intertext{whereas for $m=3k+2$}
\Pred^{(0)} &= \sum_{j=0}^{k} \sum_{i=0}^{3j+1}
q^{6j-4i+2} t^{4k+2j-2i+2} \notag \\ 
 \Pred^{(1)} &= \sum_{j=0}^{k} \sum_{i=0}^{6j}
q^{6j-2i} t^{4k+2j-2\lfloor i/2\rfloor+3}  \label{pjforktwo} \\
 \Pred^{(2)} &= \sum_{j=0}^{k-1} \sum_{i=0}^{3j+1}
q^{6j-4i+2} t^{4k+2j-2i+6}  \notag 
\end{align}
\end{conj}

\noindent
Below we summarize some checks of \eqref{predforttrim} - \eqref{pjforktwo}:

\begin{enumerate}

\item If we set $t=-1$, we recover the correct expression
for the normalized HOMFLY polynomial
\eqref{homflyfornthreea} - \eqref{homflyfornthreeb}.

\item It is easy to verify that \eqref{predforttrim} - \eqref{pjforktwo}
has the structure \eqref{torussfirst} and \eqref{torusslast},
where $S(T_{3,m})=2(m-1)$.

\item The general result \eqref{predforttrim} - \eqref{pjforktwo}
is consistent with our computations of $\Pred (T_{3,m})$ for
small values of $m$ (see examples in Sections \ref{sec-differentials}
and \ref{Sec:Examples}).

\item Taking homology with respect to $d_2$ gives the correct
result for $\khr_2 (T_{3,m})$.

\item Taking homology with respect to $d_0$ gives the correct
result for $\hfk (T_{3,m})$.

\end{enumerate}
The first three checks are fairly straightforward.  We verify the
properties (4) and (5) in the following two sections below, where we
also give the definitions of $d_2$ and $d_0$.  Another consistency
check is that $\Pred (T_{3,m})$ has the expected symmetry $\phi$.
Indeed, using the explicit form of the superpolynomial in
\eqref{predformtwo} and \eqref{predforttrim}, it is easy to verify the
following:
\begin{prop}
For $n=2$ and $n=3$, there is an involution
\begin{equation}
\label{symmphihh}
\phi: \hk_{i,j,*} (T_{n,m}) \to \hk_{i,-j,*} (T_{n,m})
\end{equation}
\end{prop}
\noindent
In other words, for torus knots $T_{2,m}$ and $T_{3,m}$,
the reduced superpolynomial $\Pred (T_{n,m})$
can be written as a polynomial in $a$, $t$,
and $y = (q^{-1} t^{-1/2} + q t^{1/2})^2$,
in agreement with the genus expansion structure \eqref{gexpansion}.


\subsection{Reduction to $\khr$}

As we explained in Section \ref{sec-differentials},
the reduction to the $sl(N)$ knot homology involves
taking cohomology with respect to the differentials $d_N$
and specializing to $a=q^N$.
Unlike the case of $(2,m)$ torus knots discussed earlier
in this section, the triply graded theory of $T_{3,m}$
is complicated enough that any differential $d_N$ can be
potentially non-zero if $m$ is sufficiently large.
In order to see this, we recall that $d_N$ is graded
of degree $(-2,2N,-1)$ for $N \geq 1$. In particular,
since it lowers the $a$-grading by 2 units and
$t$-grading by 1 unit, it should necessarily involve
the terms from $\Pred^{(1)}$ in \eqref{predforttrim}.

First, let us consider the case $m=3k+1$. 
It is convenient to split the sum over $i$ in
the expression \eqref{pjforkone} for $\Pred^{(1)}(T_{3,3k+1})$
into a sum over even and odd values of $i$, and rewrite the result as:
\begin{align}
\label{poneieven}
& \Pred^{(1)}_+ (T_{3,3k+1}) = a^{6k+2} \sum_{j=1}^{k} \sum_{i=0}^{3j-1}
q^{6j-4i-2}  t^{4k+2j-2i+1} \\
\label{poneiodd}
& \Pred^{(1)}_- (T_{3,3k+1}) = a^{6k+2} \sum_{j=1}^{k} \sum_{i=0}^{3j-2}
q^{6j-4i-4} t^{4k+2j-2i+1}
\end{align}
Now, we want to study what happens to these terms under
the action of $d_N$. Notice, here we tacitly identify
the elements of the homology groups $\hk$ with
the corresponding terms in the superpolynomial.
For example, in this terminology, a non-trivial action
of the graded differential $d_N$ is described by
a multiplication by $a^{-2} q^{2N} t^{-1}$.
Applying this to \eqref{poneieven} -- \eqref{poneiodd}
and rearranging the sum, we find
\begin{align}
\label{pzeropota}
& a^{6k} \sum_{j=N}^{k+N-1} \sum_{i=N-1}^{3j+1-2N}
q^{6j-4i} t^{4k+2j-2i} \\
\label{pzeropotb}
& a^{6k} \sum_{j=N-1}^{k+N-2} \sum_{i=N-2}^{3j+2-2N}
q^{6j-4i} t^{4k+2j-2i}
\end{align}
In this form, it is easy to recognize some of
the terms from $\Pred^{(0)} (T_{3,3k+1})$.
Indeed, comparing the range of the summation
in \eqref{poneieven} and \eqref{poneiodd}
with the one in \eqref{pjforkone}, we conclude
that $d_N$ can be potentially non-trivial
for torus knots $T_{3,3k+1}$ with $k \geq N-1$.

Similarly, we find that the terms in the expressions
\eqref{poneieven} and \eqref{poneiodd} can potentially be
in the image of $d_N$ acting on the following terms
in $\Pred^{(2)} (T_{3,3k+1})$:
\begin{align}
\label{ptwopota}
& a^{6k+4} \sum_{j=2-N}^{k+1-N} \sum_{i=2-N}^{3j+2N-2}
q^{6j-4i} t^{4k+2j-2i+4} \\
\label{ptwopotb}
& a^{6k+4} \sum_{j=1-N}^{k-N} \sum_{i=1-N}^{3j+2N-1}
q^{6j-4i} t^{4k+2j-2i+4}
\end{align}
Again, comparing these expressions with \eqref{pjforkone},
we conclude that $d_N$ has to be trivial, unless $k \geq N-1$.

Summarizing, we find that, for torus knots $T_{3,3k+1}$,
all differentials $d_N$ with $N \leq k+1$ can potentially be non-trivial.
Notice, in particular, that there are terms in
$\Pred^{(1)} (T_{3,3k+1})$ which have the right grading
to be in the image of $d_N$ as well as to map under $d_N$
to some other terms in $\Pred^{(0)}$.
Unfortunately, in this case, the structure of our triply graded theory alone
does not uniquely determine the action of $d_N$ for general $N$.
For $N=2$, we find that $d_2$ acts on the following terms
in $\Pred^{(1)}_- (T_{3,3k+1})$ and $\Pred^{(2)}(T_{3,3k+1})$:
\begin{equation}
\label{dnactsonforkone}
a^{6k+2} \sum_{j=1}^{k+2-N} \sum_{i=0}^{3j-2}
q^{6j-4i-4} t^{4k+2j-2i+1}
+ a^{6k+4} \sum_{j=0}^{k+1-N} \sum_{i=0}^{3j}
q^{6j-4i} t^{4k+2j-2i+4}
\end{equation}
and maps them to the corresponding terms in 
$\Pred^{(0)}(T_{3,3k+1})$ and $\Pred^{(1)}_+ (T_{3,3k+1})$.
Indeed, subtracting all these terms from $\Pred (T_{3,3k+1})$
and specializing to $a=q^2$, we obtain
\begin{align}
\Pred_2 (T_{3,3k+1}) &= \Pred (q^2,q,t)
- (1 + t) q^{12k} \sum_{j=1}^{k} \sum_{i=0}^{3j-2}
q^{6j-4i} t^{4k+2j-2i}   \notag \\ 
& \quad \quad \quad \quad {}- (1 + t^{-1}) q^{12k+8} \sum_{j=0}^{k-1} \sum_{i=0}^{3j}
q^{6j-4i} t^{4k+2j-2i+4} =  \label{khredtwoforkone} \\ 
&= (1 + q^4 t^2 + q^6 t^3 + q^{10} t^5)
\sum_{i=0}^{k-1} q^{6k+6i} t^{4i}
+ q^{12k} t^{4k} \notag
\end{align}
which agrees with the values of  $\khr_2(T_{3,3k+1})$ computed by
Shumakovitch and Bar-Natan. (Although we will not prove it here, this 
formula is almost certainly true in general;
 using \cite{KnotAtlas}, it can be easily checked for \(k<100\), for example.)

For $T_{3,3k+2}$, the analysis is similar. Again,
we find several possibilities for how $d_N$ might act
on various terms in the superpolynomial $\Pred (T_{3,3k+2})$:
\begin{align} \label{twotoplus} 
a^{6k+6} \sum_{j=0}^{k+1-N} \sum_{i=0}^{3j+1}
q^{6j-4i+2} t^{4k+2j-2i+6}
& ~~~\to~~~ a^{6k+4} \sum_{j=N-1}^{k} \sum_{i=N-2}^{3j-2N+2}
q^{6j-4i}  t^{4k+2j-2i+3}
\\ \label{twotominus} 
a^{6k+6} \sum_{j=0}^{k-N} \sum_{i=0}^{3j+1}
q^{6j-4i+2} t^{4k+2j-2i+6}
& ~~~\to~~~ a^{6k+4} \sum_{j=N}^{k} \sum_{i=N-1}^{3j-2N}
q^{6j-4i-2} t^{4k+2j-2i+3}
\\ \label{plustozero} 
a^{6k+4} \sum_{j=0}^{k-N+1} \sum_{i=0}^{3j}
q^{6j-4i}  t^{4k+2j-2i+3}
& ~~~\to~~~ a^{6k+2} \sum_{j=N-1}^{k} \sum_{i=N-1}^{3j-2N+2}
q^{6j-4i+2} t^{4k+2j-2i+2}
\\ \label{minustozero} 
a^{6k+4} \sum_{j=0}^{k-N+2} \sum_{i=0}^{3j-1}
q^{6j-4i-2} t^{4k+2j-2i+3}
& ~~~\to~~~ a^{6k+2} \sum_{j=N-2}^{k} \sum_{i=N-2}^{3j+3-2N}
q^{6j-4i+2} t^{4k+2j-2i+2}
\end{align}
By analogy with $(3,3k+1)$ torus knots,
one might expect that in the present case
$d_2$ acts as in \eqref{twotoplus} and \eqref{minustozero}.
In other words, one might expect that $d_2$ acts on
the following terms in $\Pred (T_{3,3k+2})$:
\begin{equation}
\label{dnactsonforktwo}
a^{6k+6} \sum_{j=0}^{k+1-N} \sum_{i=0}^{3j+1}
q^{6j-4i+2} t^{4k+2j-2i+6}
+ a^{6k+4} \sum_{j=0}^{k-N+2} \sum_{i=0}^{3j-1}
q^{6j-4i-2} t^{4k+2j-2i+3}
\end{equation}
Indeed, this 
 leads to the following result for the $sl(2)$ homology:
\begin{equation}
\label{sltworedforktwo}
\khr_2 (T_{3,3k+2}) =
(1 + q^4 t^2 + q^6 t^3 + q^{10} t^5)
\sum_{i=0}^{k} q^{6k+2+6i} t^{4i}
- q^{12(k+1)} t^{4k+5}
\end{equation}
which again agrees with the calculated value. 


\begin{remark}
  As we pointed out earlier, our prediction for $\hk (T_{3,m})$ enjoys
  a symmetry \eqref{symmphihh}, which means that the superpolynomial
  $\Pred (T_{3,m})$ can be written as a polynomial in $a$, $t$, and $y
  = (q^{-1} t^{-1/2} + q t^{1/2})^2$, in agreement with the genus
  expansion structure.  What is more surprising is that $d_N$ acts in
  a way that respects this structure!  Indeed, it easy to verify that
  both expressions in \eqref{dnactsonforkone} and
  \eqref{dnactsonforktwo} can be written in terms of the variables
  $a$, $t$, and $y$.
\end{remark}


\subsection{Reduction to $\hfk$}

We find that, for all $(3,m)$ torus knots,
the differential $d_0$ acts on the same terms as $d_2$.
Indeed,
following the same steps as in \eqref{khredtwoforkone}, we obtain for $m = 3k + 1$:
\begin{align}
\Pred_0(T_{3,3k+1}) &= \Pred (a=t^{-1},q,t)
- (1+t^{-1}) \sum_{j=1}^{k} \sum_{i=0}^{3j-2}
q^{6j-4i-4} t^{-2k+2j-2i-1} \notag \\ 
& \quad \quad {} - (1+t^{-1}) \sum_{j=0}^{k-1} \sum_{i=0}^{3j}
q^{6j-4i} t^{-2k+2j-2i} = \\ 
&= t^{-2k} \Big[ 1 + \sum_{i=1}^k \big(
q^{6i} t^{2i} + q^{6i-2} t^{2i-1} + q^{-6i+2} t^{-4i+1}
+ q^{-6i} t^{-4i} \big) \Big] \notag 
\end{align}
and for $m=3k+2$:
\begin{equation}
\Pred_0 (T_{3,3k+2}) = t^{-2k-1} \Big[
(q^2 t + 1 + {1 \over q^2t}) + \sum_{i=1}^k \big(
q^{6i+2} t^{2i+1} + q^{6i} t^{2i} + q^{-6i} t^{-4i}
+ q^{-6i-2} t^{-4i-1} \big) \Big]
\end{equation}
In both cases, this agrees with the known expressions for $\hfk(T_{3,m})$.  


\subsection{Partial results for $T_{n,m}$}

Hoping to extend the above results to all torus knots $T_{n,m}$,
one would like to have a more direct way of deriving
the superpolynomial from the general expression
\eqref{homflyfortnmnew} for the HOMFLY polynomial.
For example, our expression \eqref{predformtwo} -- \eqref{ppformtwo}
for the reduced superpolynomial of $T_{2,m}$ can be obtained
directly from the general formula \eqref{homflyfortnmnew}
for the HOMFLY polynomial by inserting powers of $(-t)$
and expanding the denominator in a power series:
\begin{equation}
\label{predtwofromhomfly}
\Pred (T_{2,m}) =
(-aqt)^{m-1} \left( {1 - q^{-2} t^{-2} \over 1 - q^{-4} t^{-2}} \right)
\left[ {1+ a^2 q^{-2} t \over 1-q^{-2} t^{-2}}
+ q^{-2m} (-t)^{2-m} {a^2 + q^{-2} t^{-3} \over 1 - q^{-2} t^{-2}} \right]
\end{equation}
Notice, the two terms inside the square brackets correspond
to $\b=0$ and $\b=1$ terms in \eqref{homflyfortnmnew}.
Similarly, for $n=3$, one has three terms in \eqref{homflyfortnmnew},
which correspond to $\b=0$, $1$, and $2$.
Comparing the structure of these terms with the corresponding
terms in the superpolynomial \eqref{predforttrim}, we find that,
again, certain parts of the superpolynomial can be obtained
directly from the HOMFLY polynomial.
Namely, these are the terms which correspond to $\b=0$ and $\b=2$.
They have the form similar to the $\b=0$ and $\b=1$ terms in 
\eqref{predtwofromhomfly} and suggest that, for a general
torus knot $T_{n,m}$, certain parts of the superpolynomial
are also given by a simple modification of the terms with
$\b=0$ and $\b=n-1$ in the HOMFLY polynomial \eqref{homflyfortnmnew}.
Namely, up to an overall power of $a$, $q$, and $t$,
the contribution of the $\b=0$ term to the superpolynomial
looks like
\begin{equation}
\label{bzeroterm}
\prod_{j=1}^{n-1} {1 + a^2 q^{-2j} t \over 1 - t^{-2} q^{-2(j+1)} }
\end{equation}
and the contribution of the $\b=n-1$ term looks like
\begin{equation}
\label{largestbterm}
\prod_{j=1}^{n-1} {a^2 + q^{-2j}t^{-3} \over 1 - t^{-2} q^{-2(j+1)} }
\end{equation}
where the terms in the denominator are understood to be expanded
in a power series. We analyze the contribution of the $\b=0$ term
more carefully in the following section.


\section{Stable homology of torus knots}\label{sec-torus-2}

Although we were unable to produce a general formula for the
superpolynomial of \(T_{m,n}\), we can make a prediction about its
behavior as \(m \to \infty\). To be precise, let us define
\begin{equation}
\label{Eq:Limit}
\CP_s(T_{m,n}) =  (a^{-1}q)^{(m-1)(n-1)} \CP(T_{m,n}).
\end{equation}

\begin{figure}
\centering\includegraphics[scale=1.3]{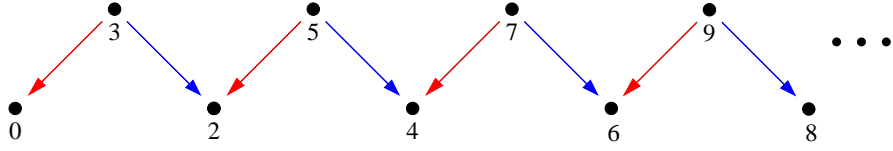}
\caption{Dot diagram for the stable
  superpolynomial of \(T_2,\) obtained as the  limit of the dot diagrams in Figure~\ref{Fig:T22k1}. }\label{Fig:2Stable}
\end{figure}

This has the effect of 
translating the dot diagram for \(\CP (T_{m,n})\) so that the leftmost
dot is always at the origin of the \((a,q)\) coordinate system.
 We then let
\begin{equation}
\CP_s(T_n)= \lim_{m \to \infty} \CP_s(T_{m,n}).
\end{equation} 
Assuming the limit exists, we refer to 
\(\CP_s(T_n)\) as the {\it stable  superpolynomial} of \(T_n\).
For example, when  \(n=2\), the calculations of Section~\ref{Subsec:T(2,m)} show that the stable 
superpolynomial is given by 
\begin{equation}
\CP_s(T_2) = (1+a^2q^2t^3) \sum_{i=0}^{\infty} q^{4i}t^{2i}.
\end{equation}
As a dot diagram, this would be represented by an up and down
chain of dots, starting at the origin of coordinates and carrying on
indefinitely to the right. This is illustrated in
 Figure~\ref{Fig:2Stable}. 

\begin{conj}
\label{Conj:StabTorus}
For all \(n\), the limit of \eqref{Eq:Limit} exists and is given by
\begin{equation}
\label{Eq:StabTorus}
\Pred_s(T_n) =
\frac{(1+a^2q^2t^3)(1+a^2q^4t^5)\ldots(1+a^2q^{2n-2}t^{2n-1})}
{(1-q^4t^2)(1-q^6t^4)\ldots(1-q^{2n}t^{2n-2})}
\end{equation}
where terms in the denominator are understood to be expanded as
a series in positive powers of \(q\) and \(t\). 
\end{conj}

Of course, we should verify that if  we substitute \(t=-1\), our
prediction for the stable superpolynomial reduces to the stable HOMFLY
polynomial of \(T_n\). 

\begin{lemma}
If $\Pred_s(T_n)$ is the expression given in \eqref{Eq:StabTorus} then
\begin{equation}
\Pred_s(T_n) \vert_{t=-1} = \lim_{m \to \infty} (qa^{-1})^{(m-1)(n-1)} P(T_{m,n}).
\end{equation}
\end{lemma}

\begin{proof}
Using the formula given in 
\eqref{homflyfortnmnew} together with the symmetry \(P_K(a,q) = P_K(a,q^{-1})\), we see that
\begin{equation}
\label{Eq:Jones2}
(qa^{-1})^{(m-1)(n-1)} P(T_{m,n}) =
\frac{1-q^2}{1-q^{2n}} \sum _{\substack{
\beta+\gamma = n-1 \\ \beta,\gamma \geq 0}} q^{2m \beta} 
\left(\prod_{i=1}^\beta \frac{a^2-q^{2i}}{1-q^{2i}}\right)
\left(\prod_{j=1}^\gamma \frac{1-a^2q^{2j}}{1-q^{2j}}\right).
\end{equation}
As \(m \to \infty\), all terms of the sum will contribute higher and
higher powers of \(q\), with the exception of the one for which
\(\beta = 0 \). We thus find
\begin{align}
P_s(T_n) & = \lim_{m \to \infty} (qa^{-1})^{(m-1)(n-1)} P(T_{m,n}) \\ 
\label{Eq:StabHOM}
&= 
\frac{1-q^2}{1-q^{2n}}  \left( \prod_{j=1}^{n-1} 
\frac{1-a^2q^{2j}}{1-q^{2j}} \right) . 
\end{align}
which agrees with the expression obtained by substituting \(t=-1\) in
\eqref{Eq:StabTorus}. 
\end{proof}

Observe that our conjectured expression for the stable superpolynomial 
has the minimum size dictated by the stable HOMFLY polynomial. Indeed,
it is easy to see from equation ~\ref{Eq:StabTorus} that the
homological grading of any term in \(\CP_s(T_n)\) is congruent 
 to half its \(a\)-grading mod \(
2\) . Thus if we substitute \(t=-1\), all
terms with a given power of \(a\) will have the same sign. In
contrast, the $sl(2)$ Khovanov homology of a torus knot is usually much larger
than the minimum size predicted by its Jones polynomial. 

\subsection{Origin of the conjecture}
Conjecture~\ref{Conj:StabTorus} was derived from the following
geometric ansatz, which is in many ways more revealing.  Here,
$\hk_s(T_n)$ is the homology groups with Poincar\'e polynomial
$\Pred_s(T_n)$.

\begin{conj}
\(\hk_s(T_n)\) is the smallest  complex satisfying the following
properties:

\begin{enumerate}
\item  \(\CP_s(T_n) \in \Z[a,q,t]\). 
\item \(\hk_s(T_n)\) contains \( \hk_s (T_{n-1})\)  as a subcomplex.
\item \(\hk_s(T_n)\) is acyclic with respect to \(d_{-1},d_{-2}, \ldots,
  d_{-n+1}\). 
\item The homology of \(\hk(T_n)\)   with respect to \(d_1\) is
  one-dimensional and generated by the monomial \(1\) appearing in
  \(\CP_s(T_n)\). 
\end{enumerate}
\end{conj}

To illustrate how \eqref{Eq:StabTorus} is derived from these
properties, consider the simplest case, when \(n=2\). We start off the stable
superpolynomial with the term \(1\), which generates the homology with
respect to \(d_1\). By property (1), \(d_{-1} (1) = 0 \).
Thus for the homology with respect to \(d_{-1}\)
to vanish, we must add a term \(a^2q^2t^3\). Next, we must kill this
new term under \(d_1\). If it is in the image of \(d_1\),
anticommutativity of \(d_1\) and \(d_{-1}\) will force \(1\) to be in
the image of \(d_{-1}\), which violates property (4). Thus we are
forced to add a third term \(q^4t^2\)  which is in the image of
\(a^2q^2t^3\) under \(d_1\). 

At this point, all the hypotheses are satisfied, 
with the exception of the fact that \(q^4t^2\) is not killed by
\(d_{-1}\). Thus we are in the same situation we started out at, only
shifted over by a factor of \(q^4t^2\). Repeating the arguments
above, we see that we must add \(a^2q^6t^5 + q^8t^4\), then
\(a^2q^{10}t^7+q^{12}t^6\), and so on indefinitely. Thus the stable
superpolynomial has the form 

\begin{align}
\CP_s (T_2) & = 1 + (a^2q^2t^3+q^4t^2)\sum_{i=0}^ \infty(q^4t^2)^i \\
          & = 1 + \frac{a^2q^2t^3+q^4t^2}{1-q^4t^2} \\
          & = \frac{1 + a^2q^2t^3}{1-q^4t^2}. 
\end{align}

\begin{figure}
\centering\includegraphics{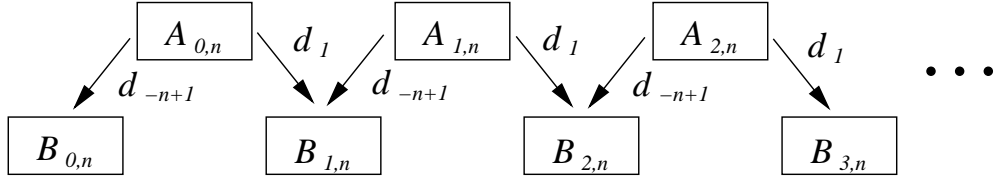}
\caption{\label{Fig:Tn} Schematic diagram of the stable complex for
  \(T_n\).  Although we've drawn each \(A_{i,n}\) and
  \(B_{i,n}\) as a finite box, they actually extend indefinitely to
  the right of the diagram.}
\end{figure}

The general case is not much different. By property (2), we start out
with \(\hk_s(T_{n-1})\), which we may inductively assume satisfies
properties (1)--(4), except that it is not acyclic with respect to
\(d_{-n+1}\). \(d_{-n+1}\) is triply graded of degree \((-2,-2n+2,-2n+1)\),
so in order to kill \(\hk(T_{n-1})\) we must add another copy of it
shifted up by \( (2,2n-2,2n-1)\). The result is acyclic with
respect to \(d_{-i}\) for \(0 < i \leq n-1\), but has the wrong
homology with respect to \(d_1\). To rectify this, we add another
copy of \(\hk_s(T_{n-1})\), shifted  by \((-2,2,-1)\) relative to
the second copy. We are now back where we started, but shifted over by 
\((0,2n,2n-2)\). Repeating, we see that \(\hk_s(T_n)\) has the general
form shown in Figure~\ref{Fig:Tn}, where the
blocks labeled \(A_{i,n}\) and \(B_{i,n}\)
 each represent an appropriately shifted 
 copy of \(\hk_s(T_{n-1})\). 

We compute
\begin{align}
\CP_s(T_n) & = \CP_s(T_{n-1}) \left(1+( a^2q^{2n-2}t^{2n-1} + q^{2n}t^{2n-2}) \sum
_{i=0}^\infty (q^{2n}t^{2n-2})^i \right) \\
         & =\CP_s(T_{n-1}) \left( 1 + \frac{a^2q^{2n-2}t^{2n-1} +
  q^{2n}t^{2n-2}}{1-q^{2n}t^{2n-2}}\right) \\ 
         & = \CP_s(T_{n-1}) \left( \frac{1+
  a^2q^{2n-2}t^{2n-1}}{1-q^{2n}t^{2n-2}}\right)
\end{align}
which clearly gives the formula of \eqref{Eq:StabTorus}. 


\subsection{Reduction to $\hfk$}

Currently, it is difficult to compute \(\khrn (T_{m,n})\) for
values of \(m,n\), and \(N\) which are all larger than \(2\), so we
have no way to check Conjecture~\ref{Conj:StabTorus}
directly. As an indirect check, however,
 we can compare the homology of \(\hk_s(T_n)\)
with respect to \(d_0\) and \(d_2\) to what is known about the knot
Floer homology and $sl(2)$ Khovanov homology of torus knots. 

The stable knot Floer homology of \(T_n\) is easily
calculated from its stable Alexander polynomial. When we substitute \(a=1\)
into the formula for the stable HOMFLY polynomial in
\eqref{Eq:StabHOM}, all the terms in the product cancel, and we are
left with
\begin{align}
\Delta_s (T_n) & = \frac{1-{q^2}}{1-q^{2n}} \\ 
&= (1-q^2)\sum_{i=0}^\infty q^{2ni}. 
\end{align}
Using Ozsv{\'a}th and Szab{\'o}'s calculation of \(\HFK\) for torus
knots in  \cite{OSlens}, it follows that
that 
\begin{equation}
\label{Eq:StabHFK}
\hfk_s (T_n) = (1+q^2t) \sum_{i=0}^\infty q^{2ni} t^{2(n-1)i}. 
\end{equation}

We want to define a differential \(d_0\) on \(\hk_s (T_n)\) which
anticommutes with the other \(d_i\)'s and whose homology is given by
the expression above. As in the construction of \( \hk_s(T_n)\), we
proceed inductively. When \(n=2\), \(d_0\) is necessarily trivial.
For general \(n\), we refer to the schematic diagram of \(\hk_s (T_n)\) in
Figure~\ref{Fig:Tn}. 
By the induction hypothesis, we can assume that we've already constructed the
  differential \(d_0\) on each block. 

To describe the part of \(d_0\) that goes between blocks, observe that 
\(\hk_s (T_{n-1})\) has a subcomplex \(C_{n-1}\) obtained by omitting
\(A_{0,n-1}\) and \(B_{0,n-1}\) from the analogous diagram for \(\hk_s
(T_{n-1})\). There is a chain map
 \(\psi : \hk_s(T_{n-1}) \to \hk_s(T_{n-1}) \) which
shifts the entire complex over one unit to the right, and which
defines an isomorphism from \( \hk_s(T_{n-1}) \) to \(C_{n-1}\). We
define the component
of \(d_0\) which maps \(A_{i,n}\) to \(B_{i,n}\) to be given by
\(\psi\), and assume that all other components of \(d_0\) between the
blocks are trivial. 

First, we should check that \(d_0\) has the correct grading. The
grading of \(A_{i,n}\) is shifted by a factor of
\((2,2n-2,2n-1)\) relative to that of
\(B_{i,n}\), while the grading of \(C_{n-1}\) is shifted by
\begin{equation}
(2,2n-4,2n-3)+(-2,2,-1) = (0,2n-2,2n-4).
\end{equation}
 Thus \(d_0\) shifts the grading by \((-2,0,-3)\), as it should. 

It follows easily from the definition that \(d_0\) anticommutes
with the other differentials. Thus it only remains to check that it
has the correct homology. To see this, note that with respect to
\(d_0\), \(\hk_s(T_n)\) decomposes as a direct sum of complexes \(D_{i,n}\),
where as a group, \(D_{i,n} = A_{i,n} \oplus B_{i,n}\). Since \(D_{i,n}\) is
just \(D_{0,n}\) shifted over by a factor of \((q^{2n}t^{2n-2})^i\),
we see that the Poincar{\'e} polynomial of the homology with respect
to \(d_0\) is 
\begin{equation}
\CP_0(T_n) = \CP(D_{0,n}) \sum_{i=0}^\infty q^{2ni} t^{2(n-1)i}.
\end{equation}
 On the other hand, it
follows from the definition of \( \psi\) that 
\(H_*(D_{i,n},d_0) \cong H_*(D_{i,n-1},d_0)\), so 
\begin{equation}
\CP(D_{0,n}) =\CP(D_{0,n-1})=\ldots= \CP(D_{0,2}) = 1+a^2 q^2 t^3.
\end{equation}
 Finally, we substitute \(a=1/t\) to
obtain
\begin{equation}
\CP_0(T_n) = (1+q^2t) \sum_{i=0}^\infty q^{2ni} t^{2(n-1)i}
\end{equation}
which agrees with \eqref{Eq:StabHFK}. 


\subsection{Reduction to $\khr_2$} 
\label{Subsec:KHR2}

As a final check on Conjecture~\ref{Conj:StabTorus}, we use it to make
some predictions about the $sl(2)$ Khovanov homology of torus knots. Although
there is not a huge amount of data with which to compare our
predictions, what there is provides some of the most convincing 
 evidence for our conjectures. Our results match perfectly with the
 known computations, which had previously seemed quite difficult to 
explain. 

To predict \(\khs(T_n)\), we must understand the action of \(d_2\) on
\(\hk_s(T_n)\). As in the previous sections, we proceed inductively,
starting with \(n=2\). In this case, \(d_2\) must vanish for
dimensional reasons, and we obtain the formula for the stable Khovanov
homology simply by substituting \(a=q^2\) and \(n=2\) into
\eqref{Eq:StabTorus}:
\begin{align}
\khs(T_2) & = (1+q^6t^3)(1+q^4t^2 + q^8t^4+\ldots) \\
& = 1 + q^4t^2 + q^6t^3 + q^8t^4 + q^{10}t^5 + \ldots
\end{align}
The Khovanov homology of \(T_{2,m}\)
 is given by 
\begin{equation}
\khr_2(T_{2,m}) = q^{m-1}(1+q^4t^2 + q^6t^3 + \ldots + q^{2m} t^m).
\end{equation}
After shifting by \(q^{1-m}\), 
this agrees with the stable homology up through terms of degree
\(q^{2m}\). In general, we expect that 
 \(q^{-mn}Kh(T_{m,n})\) should  also agree  with \(Kh_s(T_n)\) in
 degrees up to \(q^{2m}\). Indeed, if we substitute \(a=q^2\), the
 lowest degree term appearing in the expression \eqref{Eq:Jones2} for \(P(T_{n,m})\)
  which does not come from the term where \(\beta=0\) 
 is \(q^{2m+2}\). 

\begin{figure}
\centering\includegraphics{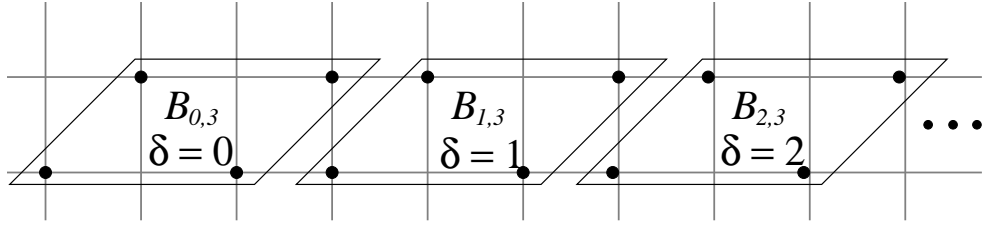}
\caption{\label{Fig:T3} What's left in \(\hk(T_3)\) after taking
  homology with respect to \(d_2\). Four generators from
  each \(B_{i,3}\) survive.}
\end{figure}

Next, we consider the case \(n=3\). Referring to Figure~\ref{Fig:Tn},
we observe that since \(d_{-2}\) lowers the \(\delta\)-grading by
\(1\) and \(d_{1}\) preserves it,
 the \(\delta\)-grading of all terms in \(A_{i,3}\) is
\(i+1\), while the \(\delta\)-grading of \(B_{i,3}\) is \(i\). Now
\(d_2\) lowers \(\delta\) by \(1\), so the only possible components of
\(d_2\) go from \(A_{i,3}\) to \(B_{i,3}\), from \(A_{i+1,3}\) to
\(A_{i,3}\), and from \(B_{i+1,3}\) to \(B_{i,3}\). In particular
\(\CF_k = \bigoplus_{i<k} D_{i,3}\) 
defines a filtration with respect to \(d_2\). We compute using the
 spectral sequence associated to this filtration. 
The differential on the \(E_0\) term is given by the restriction of
 \(d_2\) to  \(D_{i,3}\). We hypothesize that \(d_2 : A_{i,3} \to B_{i,3}\) is
nontrivial and compute its image.
Now \(A_{i,3}\) is isomorphic to \(B_{i,3}\), but shifted
in grading by \((2,4,5)\), and \(d_2\) shifts grading by
 \((-2,4,-1)\). Thus
the image of \(A_{i,3}\) under \(d_2\) will be isomorphic to
\(B_{i,3}\), but shifted by \((0,8,4)\), and the homology in
\(D_{i,3}\) will be generated by the first four terms in \(B_{i,3}\).
The Poincar{\'e} polynomial of the \(E_1\) term  is given by 
\begin{align}
\CP(E_1) &= \sum_{i=0}^\infty \CP(D_{i,3}) \\
& = \CP(D_{0,3}) \sum_{i=0}^\infty q^{6i}t^{4i} \\
\label{Eq:T3} 
& = 
\frac{(1+a^2q^2t^3 +q^4 t^2 + a^2 q^{6}t^5) }{1-q^6t^4}.
\end{align}  
This illustrated in Figure~\ref{Fig:T3}.
For dimensional reasons, there can be no further
differentials. Substituting \(a=q^2\) in \eqref{Eq:T3}, we
obtain
\begin{equation}
\label{Eq:KhT3}
\khs(T_3) = \frac{(1 +q^4 t^2  + q^6t^3 + q^{10}t^5) }{1-q^6t^4}.
\end{equation} 
This expression agrees with the pattern observed from direct
computation. For example, Figure~\ref{Fig:T(3,8)} shows
\(\khr_2(T_{3,8})\), courtesy  of Shumakovitch \cite{KhoHo}. 
As expected, the
homology agrees with \eqref{Eq:KhT3} up through powers of 
\(q^{30}\)  (here \(30 = 14 + 2\cdot8\) ).
\begin{figure}

\input T38-homology-table.tex

\vspace{-0.6cm}

\caption{The reduced Khovanov homology of $T_{3,8}$.  Here the horizontal axis corresponds to $t$, and the vertical axis to $q$.} \label{Fig:T(3,8)}
\end{figure}

In comparing these figures, it is convenient
to label generators by their
\(\delta\)-grading, since this tells us which diagonal they lie on. 
For example, the first four generators of Figure~\ref{Fig:T3} have
\(\delta\)-grading zero. They correspond to the four generators on the
highest occupied diagonal in Figure~\ref{Fig:T(3,8)}.
 The next four generators have \(\delta=1\)
and lie on the next diagonal, and so forth. 

The case \(n=4\) is somewhat more complicated. To simply things, we
assume that  as in the previous case, 
\(\CF_k = \bigoplus_{i<k} D_{i,4}\) defines a filtration
with respect to \(d_2\). Thus we are again faced with the problem of
determining the component of \(d_2\) which maps \(A_{i,4}\) to
\(B_{i,4}\). The situation is illustrated in  Figure~\ref{Fig:d2T4}. 
Possible differentials are indicated by arrows. If we assume these are
all nontrivial in rational homology (in integral homology, they are
most likely given by multiplication by \(2\)), we arrive at the following
expression for the Poincar{\'e} polynomial of \(D_{0,4}\): 
\begin{figure}
\centering\includegraphics[scale=0.9]{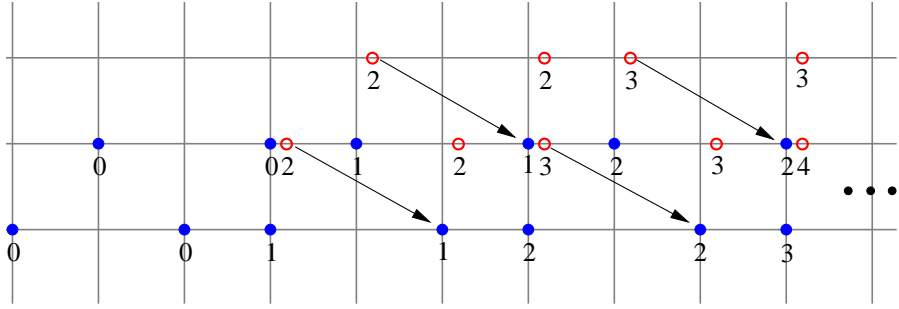}
\caption{Component of \(d_2\) from \(A_{0,4}\)
  (hollow red circles) to
  \(B_{0,4}\) (solid blue circles). Possible differentials (which we
  assume are all nonvanishing) are shown by arrows. The labels
  beneath each generator show the value of \(\delta\). }\label{Fig:d2T4} 
\end{figure}

\begin{equation}
\CP(D_{0,4})  = (1+a^2q^2t^3) \left[ 1+ q^4t^2 +
  \frac{q^6t^4(1+a^2q^4t^5)}{1-q^6t^4} \right].
\end{equation}
As before, it is easy to see there can be no further differentials, so
summing up the contributions from all \(D_{i,4}\)
we get the following prediction:
\begin{equation}
\khs(T_4) = \frac{1+q^6t^3}{1-q^8t^6} \left[ 1+ q^4t^2 +
  \frac{q^6t^4(1+q^8t^5)}{1-q^6t^4} \right].
\end{equation}
For comparison, Figure~\ref{Fig:T(4,7)} shows the Khovanov homology
of \(T_{4,7}\), again computed by \cite{KhoHo}.  
\begin{figure}

\input 47torus.tex 

\caption{The reduced Khovanov homology of $T_{4,7}$.  Here the horizontal axis correpsonds to $t$, and the vertical axis to $q$.}
\label{Fig:T(4,7)}
\end{figure}
We leave it to the reader to check that the part of the homology 
in degrees less than or equal to \( 18 + 2\cdot 7 = 32\)
  agrees with the expression above. 

As a final test, we compare with Bar-Natan's calculation of
\(\kh_2(T_{5,9})\) \cite{DBNWeb}.
 Rather than  computing a general formula for \(n=5\), we
 simply write out enough of the complex to give us the stable homology
 up to powers of \(q^{24}\). The results of the calculation are
 illustrated in Figure~\ref{Fig:StableT5}. The top half of the figure
 shows potential differentials between \(A_{0,5}\) (hollow red circles)
 and \(B_{0,5}\) (solid blue circles). Again, we assume that all these
 differentials induce nontrivial maps on rational cohomology. The
 bottom half of the figure shows the generators of \(\hk_s(T_5)\)
 which survive after taking homology. The dashed lines indicate their
 \(q\)-grading after we substitute \(a = q^2\).  By way of comparison, 
 Figure~\ref{Fig:T(5,9)} shows what we expect is 
 the reduced Khovanov homology of \(T_{5,9}\), based on Bar-Natan's
 calculation of the unreduced homology. As expected, the two agree up
 through powers of \(q^{50}\). (\(50=32+2 \cdot 9\).)

\begin{figure}
\input 59torus_delta.tex 

\vspace{-0.7cm}

\caption{\label{Fig:T(5,9)} The reduced Khovanov homology of
\(T_{5,9}\), derived from   \cite{DBNWeb}. To save space, we
 plot generators versus their \((q,\delta)\)-gradings, rather than
 \((t,q)\) as in the previous figures. }
\end{figure}

\begin{figure}

\

\vspace{2cm}
\includegraphics{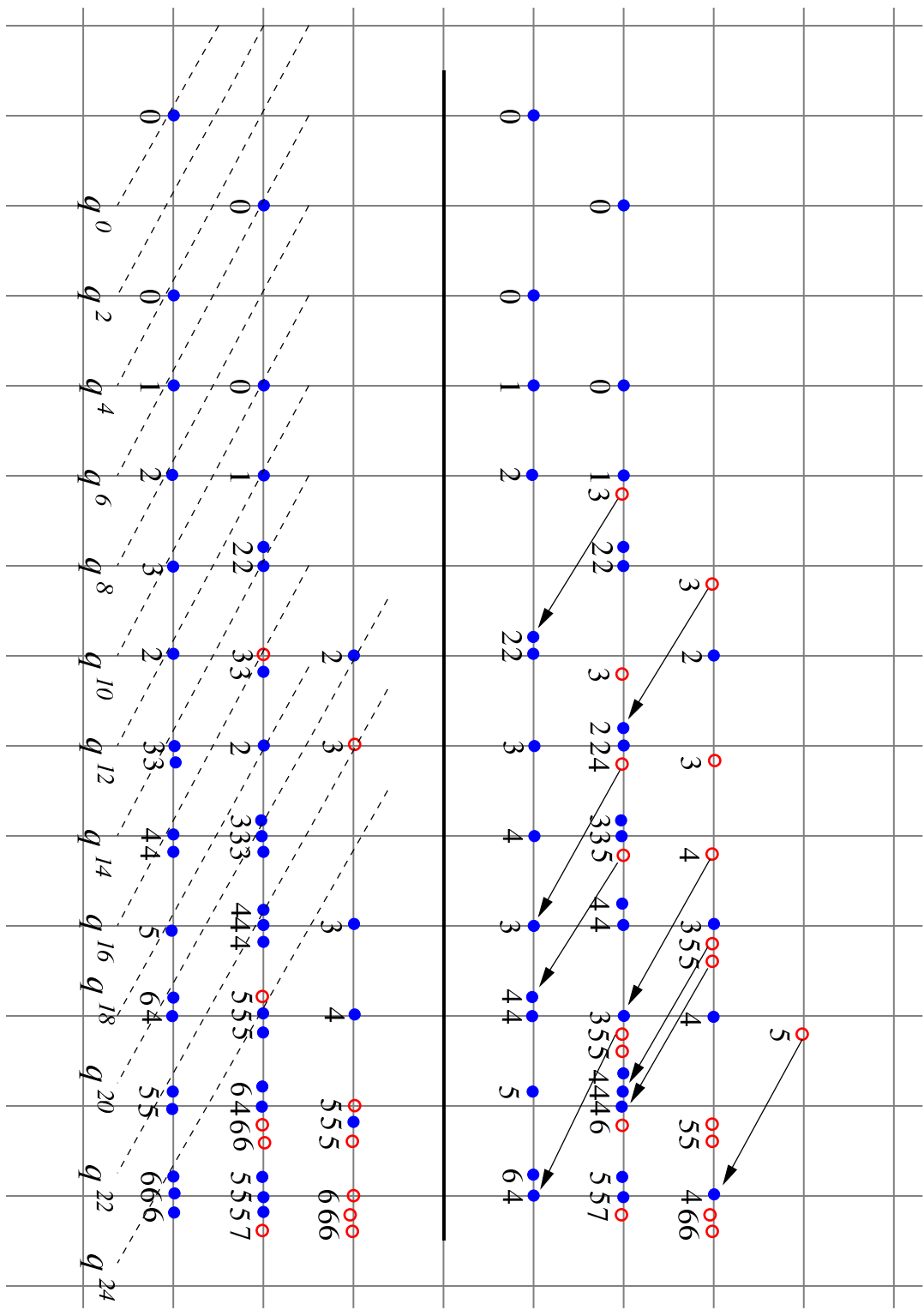}

\caption{\label{Fig:StableT5}
Computing the stable Khovanov homology of \(T_5\). Generators are
labeled by their \(\delta\)-grading}
\end{figure}

\ifarxiv
\clearpage

\section{Dot diagrams for 10 crossing knots}\label{sec-dot-diagrams}

\vspace{0.3cm}

\dotdia{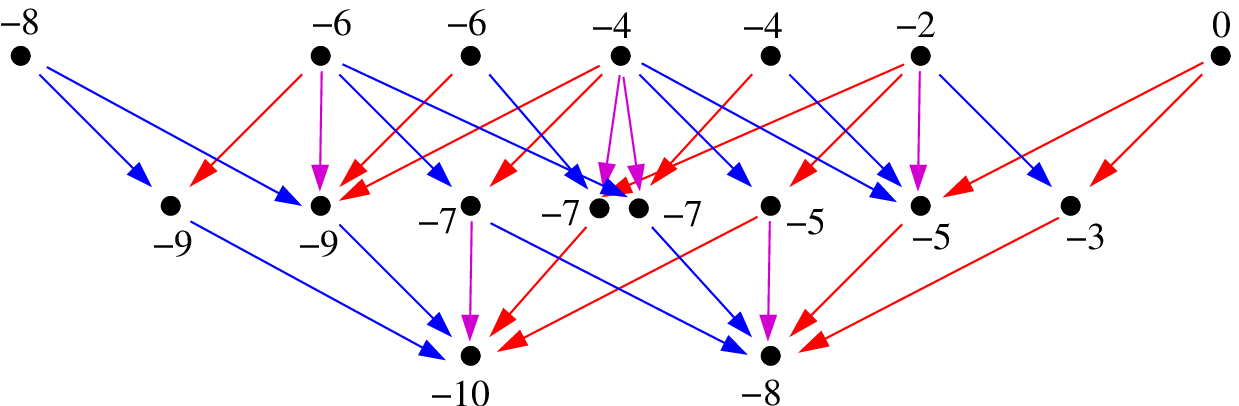}{10_{124} = \mathrm{mirror}(T_{3,5})}{-12}{1.1}

\vfill

\dotdianod{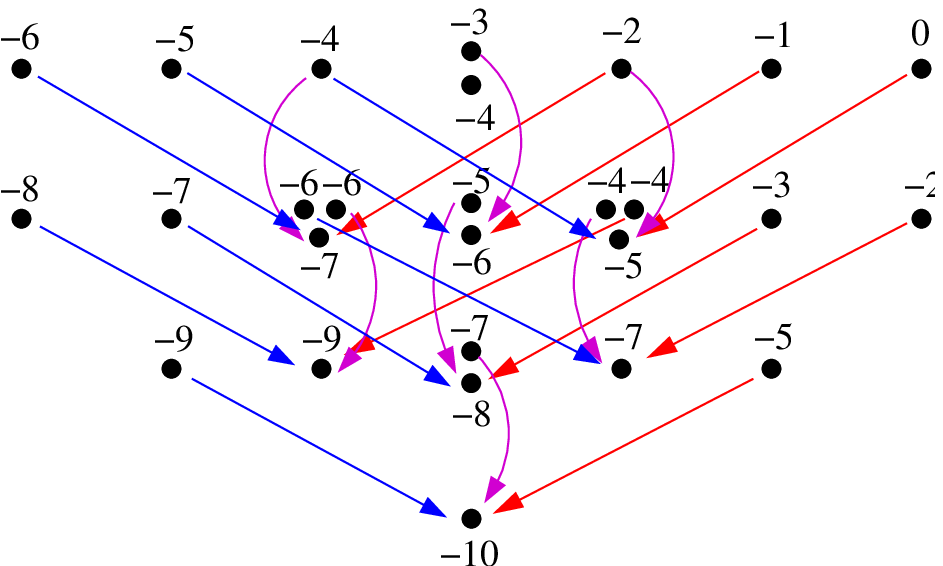}{10_{128}}{-12}{1.1}

\vfill 

\dotdia{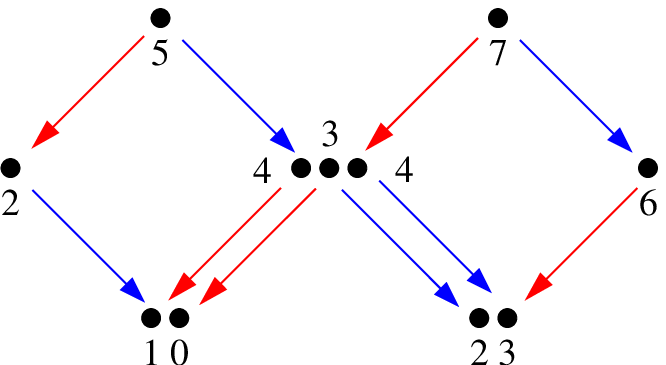}{10_{132}}{2}{1.1}

\pagebreak

\dotdia{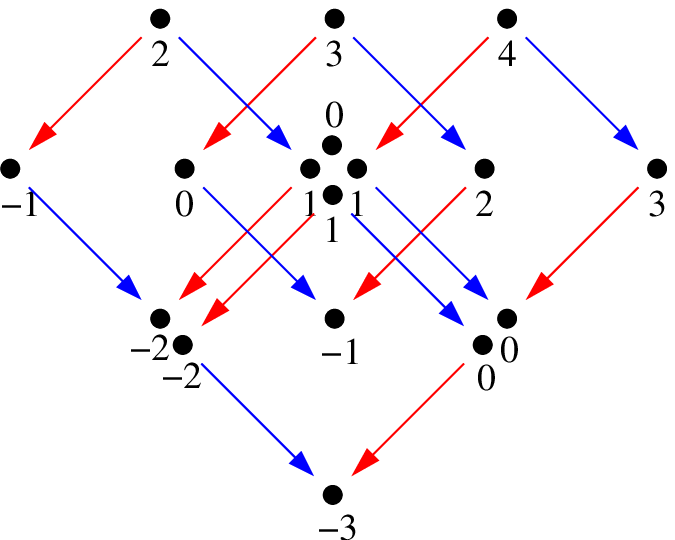}{10_{136}}{-4}{1.1}

\vfill 

\dotdia{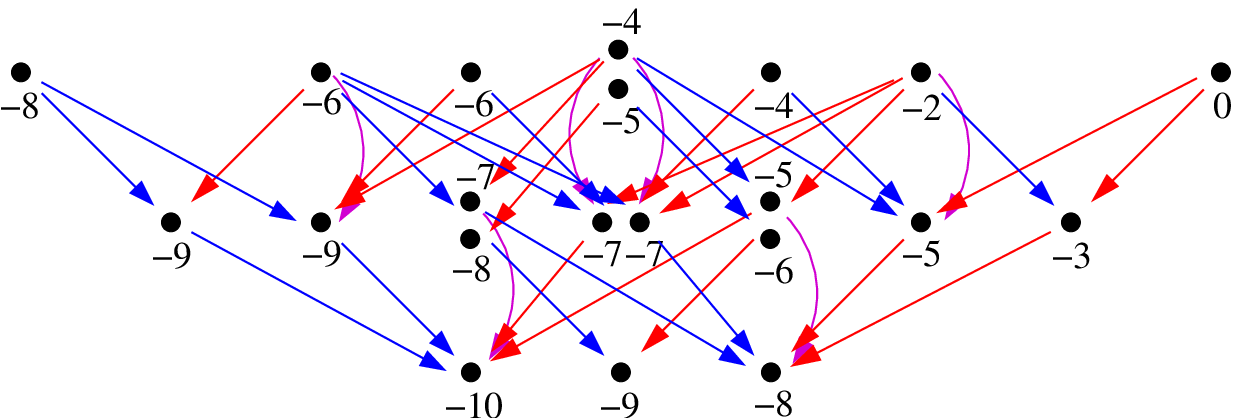}{10_{139}}{-12}{1.1}

\vfill 

\dotdia{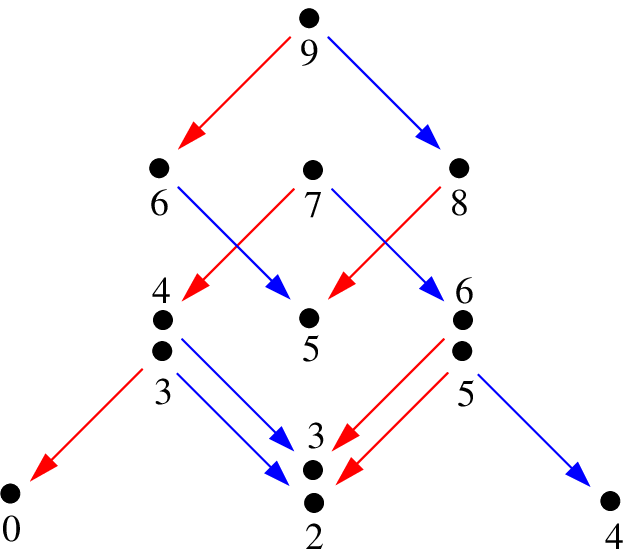}{10_{145}}{4}{1.1}

\pagebreak

\

\vspace{0.3cm}

\dotdianod{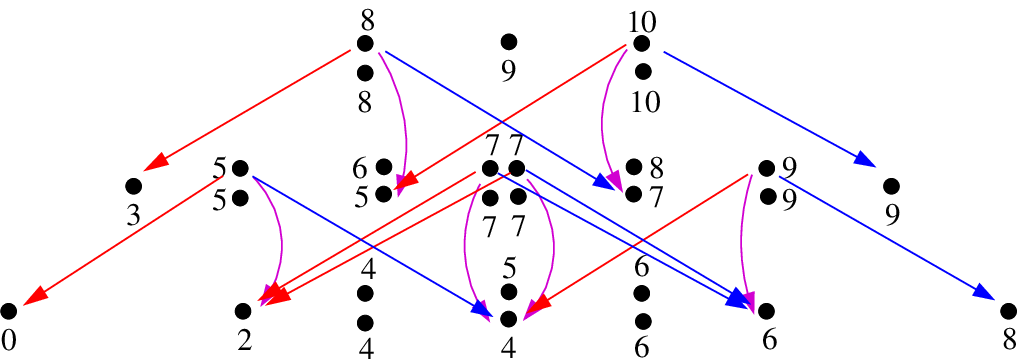}{10_{152}}{8}{1.6}

\vfill

\dotdia{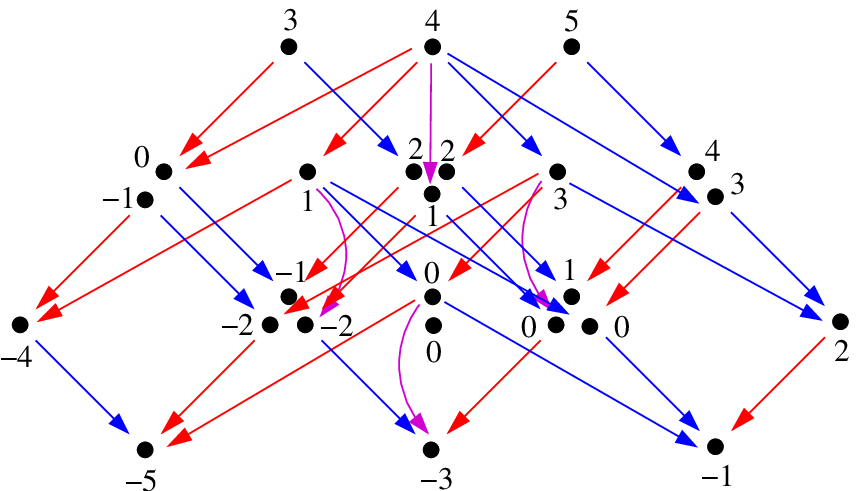}{10_{153}}{-2}{1.6}

\vspace{0.3cm}

\pagebreak

\

\vspace{0.3cm}

\dotdia{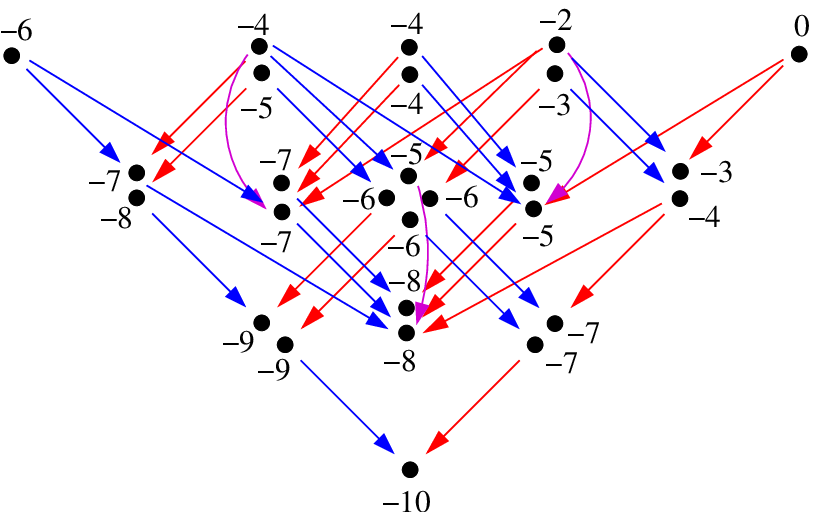}{10_{154}}{-12}{1.6}

\vfill

\dotdia{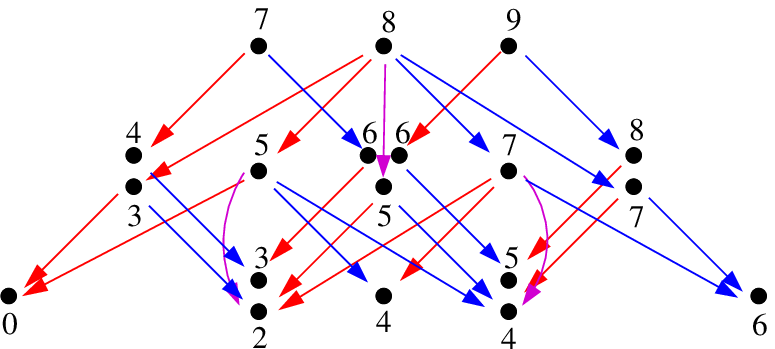}{10_{161}}{6}{1.6}

\vspace{0.3cm}
\pagebreak

\else

\section{Dot diagrams for 10 crossing knots}\label{sec-dot-diagrams}

\dotdia{10_124.eps}{10_{124} = \mathrm{mirror}(T_{3,5})}{-12}{1.1}

\dotdianod{10_128.eps}{10_{128}}{-12}{1.1}

\dotdia{10_132.eps}{10_{132}}{2}{1.1}

\dotdia{10_136.eps}{10_{136}}{-4}{1.1}

\dotdia{10_139.eps}{10_{139}}{-12}{1.1}

\dotdia{10_145.eps}{10_{145}}{4}{1.1}

\dotdianod{10_152.eps}{10_{152}}{8}{1.6}

\dotdia{10_153.eps}{10_{153}}{-2}{1.6}

\vspace{0.3cm}

\dotdia{10_154.eps}{10_{154}}{-12}{1.6}

\dotdia{10_161.eps}{10_{161}}{6}{1.6}

\fi


\ifarxiv
\addtolength{\topmargin}{-0.5in}
\addtolength{\footskip}{0.6in}
\enlargethispage{0.7in}
\fi

\end{document}